\documentclass[final,3p,times]{elsarticle}
\usepackage{amssymb}
\usepackage{graphicx,setspace}
\usepackage{booktabs}
\usepackage{psfrag}
\usepackage{amsfonts}
\expandafter\let\csname equation*\endcsname\relax
\expandafter\let\csname endequation*\endcsname\relax
\usepackage{amsmath}
\usepackage{amssymb, amsthm}
\usepackage{mathrsfs}
\usepackage{epstopdf}
\usepackage{float}
\usepackage{lineno,hyperref}
\usepackage{color}
\usepackage{amsmath}
\usepackage{dsfont} 
\usepackage{amssymb} 
\usepackage{graphicx,color}
\usepackage{extarrows}
\usepackage{xcolor}
\usepackage{subcaption}

\usepackage{graphicx}
\usepackage{amsmath,amssymb,amsthm,mathtools}
\usepackage{float}
\usepackage{booktabs}
\usepackage{hyperref}
\usepackage{cleveref}
\newcommand{\R}{\mathbb{R}}
\newcommand{\E}{\mathbb{E}}
\newcommand{\Prob}{\mathbb{P}}

\newtheorem{proposition}{Proposition}[section]

\journal{arXiv}

\begin{document}

\begin{frontmatter}

\title{
State-Dependent Saving and Heavy-Tailed Shocks in a Stochastic
Solow Growth Model via a Logistic L\'evy Framework
}

\author{Wen Wang\fnref{addr1,addr2}}\ead{13165489955@163.com}

\author{Shenglan Yuan\corref{cor1}\fnref{addr1}}\ead{shenglanyuan@gbu.edu.cn}\cortext[cor1]{Corresponding author}

\address[addr1]{\rm Department of Mathematics, School of Sciences, Great Bay University, Dongguan 523000, China }
\address[addr2]{\rm School of Mathematics and Statistics, The University of New South Wales, Sydney NSW 2052,
Australia}

\begin{abstract}
This work develops a stochastic Solow growth model in which the aggregate
saving rate is state dependent and evolves on a faster time scale than
capital. The saving target follows a mean-reverting process driven by
symmetric $\alpha$-stable L\'evy noise, while the saving control adjusts
toward this target and capital accumulates according to a Solow equation.
After scaling, the model is a three-dimensional stochastic dynamical system
with multiplicative heavy-tailed noise and a logistic saving function. We
establish invariance of the economically meaningful state space for the
discrete reflected model, derive the stationary distribution of the saving
target in the Gaussian limiting case, and characterize how the tail index
and logistic sensitivity affect the coupled dynamics. Monte Carlo
simulations show bimodal capital and saving distributions, asymmetric
first-passage behavior between low- and high-capital regimes, and a
nonlinear relationship between output volatility, the logistic sensitivity,
and the tail index. In particular, heavier tails increase the dispersion of
the scaled saving target and low-to-high first-passage frequencies, but do
not necessarily increase relative output volatility. The results highlight
the importance of state-dependent saving and non-Gaussian shocks for growth
dynamics.
\end{abstract}


\end{frontmatter}

\section{Introduction}

The saving rate is one of the two proximate determinants of long-run
economic growth, the other being the efficiency with which saved resources
are allocated to productive investment. Since Solow's seminal contribution
\cite{solow1956}, the overwhelming majority of neoclassical growth models
have treated the aggregate saving rate as a constant parameter---a fixed
fraction of income that does not respond to the state of the economy. This
modeling choice has been enormously productive analytically, but it sits
uneasily with a large body of empirical and theoretical work. Cross-country
evidence documents substantial and systematic variation in saving rates
across income levels, and the assumption of a constant saving propensity has
been shown to preclude a range of empirically relevant phenomena, including
multiple steady states and the twin-peaks pattern of income distribution.

The endogeneity of saving behavior has been a central theme in economics at
least since Kaldor's nonlinear business-cycle model \cite{kaldor1940}, in
which both investment and saving are specified as nonlinear functions of
income and capital. Kaldor's insight was that the marginal propensity to
save varies with the level of economic activity, being low at depressed
income levels and high during expansions. Modern treatments of endogenous
saving have proceeded along different lines. The optimal-growth tradition
derives saving from intertemporal utility maximization, yielding a saving
rate that depends on the capital stock and on preference parameters.
Buiter's survey of saving and endogenous growth \cite{buiter1991} emphasizes
that the determinants of the saving rate are far richer than the
constant-$s$ benchmark admits, and that different specifications of saving
behavior carry sharply different implications for long-run growth. A
parallel literature has studied state-dependent preferences and risk
aversion, showing that saving rates respond asymmetrically to changes in
wealth and risk exposure. Yet the bulk of this work remains couched either
in the representative-agent optimal-growth framework or in reduced-form
specifications that do not embed saving dynamics within a fully specified
stochastic growth process.

A more recent strand of research has begun to relax the constant-saving
assumption within the Solow framework itself \cite{Brannan2019}. A logistic saving
function---one that rises gradually from a low to a high value as income
increases---can generate multiple locally stable steady states and a
poverty-trap structure that the constant-saving model cannot produce. The
logistic form is natural because it captures the empirically observed
S-shaped relationship between saving rates and income, and because it
emerges from microfoundations with heterogeneous agents facing
non-homothetic preferences or borrowing constraints. The resulting dynamics
exhibit thresholds, hysteresis, and sensitivity to initial conditions that
are absent in the linear-saving case. However, these extensions have largely
been developed in deterministic settings. The interaction between
state-dependent saving and stochastic shocks to the growth process remains
largely unexplored.

At the same time, the stochastic growth literature has moved beyond the
Gaussian framework that dominated early extensions of the Solow model. The
standard approach replaces deterministic productivity or depreciation with a
Brownian motion, yielding a diffusion process for capital. Yuan et al.\
\cite{yuan2026} analyze the deterministic skeleton of the model,
characterizing the steady states and their stability properties as functions
of the saving-function parameters. While analytically tractable, this
specification imposes continuous sample paths and finite variance---properties
that are difficult to reconcile with the observed frequency of large,
discontinuous changes in macroeconomic aggregates. L\'evy-driven growth
models have been proposed as a more realistic alternative. Abebe et al.\
\cite{abebe2026} integrate $\alpha$-stable L\'evy noise into the Solow model
and show that heavy-tailed, jump-diffusion dynamics produce qualitatively
different long-run behavior from Gaussian perturbations, including stochastic
bifurcations and altered forecast distributions. Related work \cite{imkeller2006exit} applies
multiplicative $\alpha$-stable noise to logistic growth models in ecology and
economics, where it generates escape probabilities and mean exit times that
differ markedly from the Brownian case. These developments suggest that the
distributional characteristics of the noise process are not a technical
detail but a substantive modeling choice with first-order implications for
economic dynamics.

This work brings these strands together by developing a stochastic growth
model in which the saving rate is an endogenous, state-dependent control
variable that evolves on a faster time scale than the capital stock. The
central modeling assumption is that aggregate saving behavior responds to
the level of capital per effective worker through a logistic function, with
the sensitivity of that response controlled by a parameter $\gamma$. This
captures the idea that households and firms adjust their saving propensities
gradually as the economy moves away from a reference capital level---a
reference level that is itself determined by the average of the low and high
saving rates. The separation of time scales is not merely a technical
convenience: it reflects the observation that saving decisions, operating
through financial intermediation and portfolio adjustment, adjust more
rapidly than the physical capital stock, which is constrained by
installation and depreciation lags.

The model is formulated as a three-dimensional system in which capital per
effective labor $k(t)$ evolves slowly, while a scaled saving control $s(t)$
and its stochastic target $p(t)$ evolve on a faster scale. The target $p(t)$
follows a mean-reverting process driven by an $\alpha$-stable L\'evy noise,
so that the saving environment is subject to heavy-tailed discontinuous
shocks. The discrete-time approximation of this system, obtained by replacing
derivatives with difference quotients with explicit step $\Delta t$, yields a
stochastic difference equation system in which the saving control in period
$n+1$ responds to the stochastic target in the same period, and capital in
period $n+1$ depends on the saving control that has already adjusted. The
parameter vector is $v = (\theta, s_1, s_2, \mu, \gamma, \beta, a, \sigma,
\alpha)$, with the scaling factor $\eta = 2/(s_1 + s_2)$ determined by the
saving bounds rather than independently specified.

The paper is organized as follows. Section~\ref{sec:model} develops the
mathematical model in both dimensional and scaled forms, derives the discrete
approximation with explicit time step, and discusses the time-scale
separation. Section~\ref{sec:stochastic} incorporates the L\'evy noise and
studies the stochastic dynamics, including invariance of the state space, the
stationary distribution of the scaled saving target, the generator of the
L\'evy-driven system, and the qualitative role of the tail index $\alpha$.
Section~\ref{sec:numerical} presents numerical simulations and examines the
model's implications for the relationship between state-dependent saving and
macroeconomic volatility. Section~\ref{sec:conclusion} concludes.

\section{State-Dependent Saving and Mathematical Model}\label{sec:model}

The central modeling assumption of this work is that the aggregate saving
rate depends nonlinearly on the state of the economy and evolves on a faster
time scale than the capital stock. Let $k(t)>0$ denote capital per unit of
effective labor and let $s(t)$ be a scaled saving control variable. We model
the state-dependent saving rate by the logistic function
\begin{equation*}
	g(k;\gamma)
	=
	s_1
	+
	\frac{s_2 - s_1}{1 + e^{-\gamma (k-1)}},
\end{equation*}
where $0 \le s_1 < s_2 \le 1$ and $\gamma > 0$ controls the sensitivity of
saving behavior to deviations of capital from its reference level $k=1$. The
function $g$ is bounded between $s_1$ and $s_2$; it represents the scaled
saving \emph{target}, not the realized saving rate. The realized (scaled)
saving rate $s(t)$ can deviate from $g(k)$ because of the stochastic target
and the finite adjustment speed, as specified below.

Consider the following dimensional system:
\begin{align*}
	d\tilde{k}(t) &= \big[\tilde{s}(t)\,b\,\tilde{k}(t)^{\theta} - \mu\,\tilde{k}(t)\big]\, dt,
	\\
	d\tilde{s}(t) &= -\beta \big[\tilde{s}(t)- \tilde{p}(t)\big]\, dt,  \\
	d\tilde{p}(t) &= -a\big[\tilde{p}(t) - \tilde{g}(\tilde{k}(t))\big]\, dt
	\;+\; \sigma\,\sqrt{\tilde{p}(t)\big(1 - \tilde{p}(t)\big)}\, dL^{\alpha}(t),
\end{align*}
where
\begin{equation*}
	\tilde{g}(\tilde{k}) = s_1 + \frac{s_2 - s_1}{1 + e^{-\tilde{\gamma}(\tilde{k}-\tilde{k}^{\ast})}},
	\qquad
	\tilde{k}^{\ast} = \left( \frac{b(s_1 + s_2)}{2 \mu} \right)^{1/(1-\theta)}.
\end{equation*}

\subsection{Scaled Variables and Time-Scale Separation}

Introducing the dimensionless state variables $k(t)$, $s(t)$, and $p(t)$
using the scalings
\begin{equation*}
	\tilde{k}(t) = \tilde{k}^{\ast}\,k(t),
	\qquad
	\tilde{s}(t) = \frac{s_1+s_2}{2}\,s(t),
	\qquad
	\tilde{p}(t) = \frac{s_1+s_2}{2}\,p(t),
\end{equation*}
yields the scaled system
\begin{align}
	dk(t) &= \mu\big[ s(t)\, k(t)^{\theta} - k(t) \big]\, dt,
	\label{scaledk} \\
	ds(t) &= -\beta \big[s(t)-p(t)\big] dt, \label{scaleds} \\
	dp(t) &= -a \big[p(t)- \eta g(k(t))\big] dt + \sigma\sqrt{p(t)(\eta-p(t))}\,dL^\alpha(t),
	\label{scaledp}
\end{align}
where
\begin{equation*}
	\eta = \frac{2}{s_1+s_2}, \qquad
	\gamma = \tilde{\gamma}\,\tilde{k}^{\ast},
\end{equation*}
and the scaled target is
\begin{equation}\label{eq:SigmoidFunction}
	g(k(t)) = s_1 + \frac{s_2-s_1}{1+ e^{-\gamma(k(t)-1)}}.
\end{equation}
The factor $\eta$ appears in the nondimensional system because the original
target $\tilde p$ is bounded between $\tilde s_1 = (s_1+s_2)s_1/2$ and
$\tilde s_2 = (s_1+s_2)s_2/2$, and the scaling by $(s_1+s_2)/2$ maps these
to $\eta s_1$ and $\eta s_2$ (not to $0$ and $\eta$). The choice of
reflection bounds $0$ and $\eta$ used in the numerical scheme is therefore a
convenient enlargement of the economically relevant interval
$[\eta s_1, \eta s_2]$, not an identification of the realized saving rate
with the target bounds $s_1$ and $s_2$. All figures and tables in
Section~\ref{sec:numerical} report the scaled variables $s$ and $p$; the
realized scaled saving rate is $\tilde s = s/\eta \in [0,1]$ (in original
units), and the realized saving rate is bounded between $0$ and $1$ in
original units, not necessarily between $s_1$ and $s_2$.

Note that the parameters $\mu$, $\beta$, and $a$ have units of
$\text{time}^{-1}$, whereas the noise scale parameter $\sigma$ has units
$\text{time}^{-1/\alpha}$ under the standard symmetric $\alpha$-stable
convention $L^\alpha(t)\sim S_\alpha(t,0,1,0)$.

The logistic function $g(k)$ describes an endogenous saving rate
bounded between $s_1$ and $s_2$, with responsiveness controlled by
$\gamma$. The centering point $k^{\ast}$ is chosen as the positive
equilibrium of the model with constant saving rate
$\bar{s}=(s_1+s_2)/2$:
\begin{equation*}
	\bar{s}\,b\,(k^{\ast})^\theta = \mu k^{\ast},
	\qquad
	k^{\ast} =
	\left(\frac{b(s_1+s_2)}{2\mu}\right)^{\!\frac{1}{1-\theta}}.
\end{equation*}

\subsection{A Discrete Approximation with Explicit Time Step}\label{sec:discrete}

We obtain a discrete approximation to system \eqref{scaledk}--\eqref{scaledp}
with explicit time step $\Delta t>0$. Let $n$ denote the numerical step index,
$t_n=n\Delta t$, and write $k_n = k(t_n)$, $s_n=s(t_n)$, $p_n=p(t_n)$.
Replacing derivatives by forward difference quotients and using the
scaling of the $\alpha$-stable noise, the explicit Euler-type scheme with
reflector on $p$ reads
\begin{align}
	p_{n+1} &= \mathcal R_{[0,\eta]}\Big[
	p_n - a\Delta t\big(p_n-\eta g(k_n)\big)
	+ \sigma (\Delta t)^{1/\alpha}\sqrt{p_n(\eta-p_n)}\,\varepsilon_{n+1}
	\Big], \label{eq:disc_p}\\
	s_{n+1} &= s_n + \beta\Delta t\big(p_{n+1}-s_n\big), \label{eq:disc_s}\\
	k_{n+1} &= k_n + \mu\Delta t\big(s_{n+1}k_n^\theta-k_n\big), \label{eq:disc_k}
\end{align}
where
\begin{equation*}
	\mathcal R_{[0,\eta]}(x) = \eta - \left| (x \bmod 2\eta) - \eta \right|,
	\qquad
	\{\varepsilon_n\} \stackrel{\text{i.i.d.}}{\sim} S_\alpha(0,1,0),
\end{equation*}
and $S_\alpha(0,1,0)$ denotes the standard symmetric $\alpha$-stable
distribution with characteristic function
$\E[e^{iu\varepsilon_n}]=\exp(-|u|^{\alpha})$. The reflector
$\mathcal R_{[0,\eta]}$ is a specular (mirror) reflection. The
term $(\Delta t)^{1/\alpha}$ is the correct temporal scaling for
$\alpha$-stable increments: if $L^\alpha(t)$ is the underlying L\'evy
process, then $L^\alpha((n+1)\Delta t)-L^\alpha(n\Delta t)
\stackrel{d}{=}(\Delta t)^{1/\alpha}\varepsilon_{n+1}$.

Numerical step $n$ corresponds to time $t_n=n\Delta t$ quarters, so one
quarter comprises $1/\Delta t$ numerical steps. In the numerical experiments
of Section~\ref{sec:numerical} we take $\Delta t=0.05$, so each quarter
contains $20$ numerical updates. The parameter vector is
\begin{equation*}
	v = \left(\theta,\, s_1,\, s_2,\, \mu,\, \gamma,\, \beta,\, a,\, \sigma,\, \alpha\right),
\end{equation*}
with $\eta=2/(s_1+s_2)$ determined by the saving bounds and not an
independent parameter. It is convenient to write
$\mathbf y_n = (k_n, s_n, p_n)^T$ for the state of the system.
Stochastic bifurcations in two-time-scale dynamical systems driven by
$\alpha$-stable L\'evy noise have been investigated in \cite{yuan2021timescale,yuanwang2023,tesfay2021}.

The \emph{continuous-time} system \eqref{scaledk}--\eqref{scaledp}
with $\alpha$-stable noise and reflection at $0$ and $\eta$ should be
regarded as a formal modeling device. All rigorous statements and all
numerical experiments in this paper refer to the discrete reflected scheme
\eqref{eq:disc_p}--\eqref{eq:disc_k}.

\section{Stochastic Dynamics under L\'evy Noise}\label{sec:stochastic}

Having established the deterministic skeleton and its discrete
approximation, we now incorporate the $\alpha$-stable L\'evy noise and
analyze the resulting stochastic dynamics. Randomness enters directly only
through the scaled saving target $p$, and is transmitted to the scaled
saving control $s$ and the capital stock $k$ through the chain
$p\to s\to k$.

\subsection{The $\alpha$-Stable L\'evy Process}\label{sec:levy}

Let $L^{\alpha}(t)$ denote a standard symmetric $\alpha$-stable L\'evy
process with tail index $\alpha \in (0,2]$, scale parameter $1$, skewness
$0$, and location $0$. The process has stationary independent increments,
and for $\alpha < 2$ its sample paths are almost surely discontinuous,
exhibiting jumps of all sizes. The characteristic function of $L^{\alpha}(t)$
is
\begin{equation*}
	\E\left[ e^{i u L^{\alpha}(t)} \right] = \exp\left( -t |u|^{\alpha} \right), \qquad u \in \R.
\end{equation*}

For $\alpha=2$, the process satisfies $L_t^2\overset{d}{=}\sqrt{2}W_t$,
where $W_t$ is a standard Brownian motion. For $\alpha < 2$, the tails of
the distribution decay as a power law:
$\Prob(|L^{\alpha}(t)| > x) \sim C_{\alpha} t x^{-\alpha}$ as $x \to \infty$,
where $C_{\alpha}$ is a constant depending on $\alpha$. This heavy-tailed
property implies that the variance of $L^{\alpha}(t)$ is infinite for
$\alpha < 2$, and the mean absolute deviation is not finite for
$\alpha \le 1$.

In the discrete approximation \eqref{eq:disc_p}, the noise terms
$\varepsilon_n$ are independent and identically distributed standard
symmetric $\alpha$-stable random variables with characteristic function
$\E[e^{iu\varepsilon_n}] = \exp(-|u|^{\alpha})$. The discrete noise process
therefore inherits the heavy-tailed jump-diffusion character of its
continuous-time counterpart.

\subsection{Boundedness and Positive Invariance}\label{sec:invariance}

A preliminary concern is whether the state variables remain within
economically meaningful bounds. The scaled saving control $s_n$ and the
scaled saving target $p_n$ are designed to lie in $[0,\eta]$, and capital
$k_n$ must remain positive. The multiplicative noise structure in
\eqref{eq:disc_p} is designed to respect these bounds: the factor
$\sqrt{p_n(\eta-p_n)}$ vanishes at $p_n=0$ and $p_n=\eta$. The reflector
$\mathcal R_{[0,\eta]}$ then guarantees $p_{n+1}\in[0,\eta]$ even after a
large jump.

\begin{proposition}[Positive Invariance of the Discrete Scheme]\label{prop:invariance}
	Assume $p_0\in[0,\eta]$, $s_0\in[0,\eta]$, $k_0>0$, $\theta\in(0,1)$,
	$0<\beta\Delta t\le 1$, and $0<\mu\Delta t<1$. Then the discrete scheme
	\eqref{eq:disc_p}--\eqref{eq:disc_k} satisfies $p_n\in[0,\eta]$,
	$s_n\in[0,\eta]$, and $k_n>0$ for all $n\ge 0$ almost surely.
\end{proposition}

\begin{proof}
	The reflector gives $p_{n+1}=\mathcal R_{[0,\eta]}(\cdot)\in[0,\eta]$
	whenever the argument is finite; the multiplicative factor
	$\sqrt{p_n(\eta-p_n)}$ is well defined and the additive drift
	$-a\Delta t(p_n-\eta g(k_n))$ is finite for $k_n>0$. For $s_{n+1}$,
	the update is a convex combination
	$s_{n+1}=(1-\beta\Delta t)s_n+\beta\Delta t\,p_{n+1}$, so
	$s_{n+1}\in[0,\eta]$. For $k_{n+1}$, we have
	$$
	k_{n+1}=k_n(1-\mu\Delta t)+\mu\Delta t\,s_{n+1}k_n^{\theta}.
	$$
	Because $k_n>0$, $s_{n+1}\ge0$, and $1-\mu\Delta t>0$, the right-hand
	side is strictly positive. Induction gives the result.
\end{proof}

Proposition~\ref{prop:invariance} ensures that the discrete model as in \eqref{eq:disc_p}--\eqref{eq:disc_k} is
well posed and that the state space
$\mathcal D = (0,\infty)\times[0,\eta]\times[0,\eta]$ is positively
invariant. We do not claim a corresponding invariance result for the
continuous-time reflected L\'evy SDE; the discrete reflected scheme is the
operative model in this work.

\subsection{Stationary Distribution of the Scaled Saving Target}\label{sec:stationary}

The scaled saving target $p_n$ evolves autonomously once the path of $k_n$
is known, but in the full system \eqref{eq:disc_p}--\eqref{eq:disc_k} the two are coupled. To fix intuition, it
is useful to first analyze the behavior of $p$ for a fixed capital level
$k$, treating $g(k)$ as a constant. This provides a benchmark for the full
coupled system.

Consider the one-dimensional stochastic differential equation
\begin{equation}\label{eq:p_only}
	dp(t) = -a\left[ p(t) - \eta g \right] dt
	+ \sigma \sqrt{p(t)(\eta - p(t))} \, dL^{\alpha}(t),
\end{equation}
where $g \in (0,1)$ is a constant. This is a L\'evy-driven diffusion on the
interval $[0,\eta]$ with mean-reverting drift and state-dependent noise.

For the Gaussian case ($\alpha=2$), the stationary distribution of
\eqref{eq:p_only} can be derived explicitly. The backward generator of the
diffusion is
\begin{equation*}
	\mathcal{L} \phi(p) = -a(p - \eta g) \phi'(p) + \sigma^2 p(\eta - p) \phi''(p),
\end{equation*}
where we used $L^2(t)\overset{d}{=}\sqrt{2}W_t$. The operator $\mathcal L$
acts on sufficiently smooth test functions $\phi$. The corresponding
forward (Fokker--Planck) operator is
\begin{equation*}
	\mathcal L^* \rho = -\frac{\partial}{\partial p}[a(\eta g-p)\rho]
	+ \sigma^2\frac{\partial^2}{\partial p^2}[p(\eta-p)\rho].
\end{equation*}
If the process is ergodic, the time-dependent density $\rho(p,t)$ converges
to a stationary density $\pi(p)$ satisfying $\mathcal L^*\pi=0$.

Solving $\mathcal L^*\pi=0$ yields
\begin{equation}\label{eq:stationary_gaussian}
	\pi(p) \propto p^{\frac{ag}{\sigma^2} - 1} (\eta - p)^{\frac{a(1-g)}{\sigma^2} - 1},
\end{equation}
which is a Beta distribution on $[0,\eta]$ with parameters
\begin{equation*}
	A = \frac{ag}{\sigma^2}, \qquad B = \frac{a(1-g)}{\sigma^2}.
\end{equation*}
The stationary mean is $\eta g$, and the stationary variance is
\begin{equation*}
	\operatorname{Var}(p) = \frac{\eta^2 \sigma^2 g(1-g)}{a + \sigma^2}.
\end{equation*}
The density \eqref{eq:stationary_gaussian} is integrable provided
$A>0$ and $B>0$; it is bounded at the boundaries when $A\ge 1$ and $B\ge 1$.

For $\alpha < 2$, the situation is more complex. The generator of the
L\'evy-driven diffusion is
\begin{equation*}
	\mathcal{L} \phi(p) = -a(p - \eta g) \phi'(p)
	+ \int_{\R \setminus \{0\}} \left[
	\phi(p + h(p) y) - \phi(p) - h(p) y \phi'(p) \mathbf{1}_{|y| < 1}
	\right] \nu(dy),
\end{equation*}
where $h(p)=\sigma\sqrt{p(\eta-p)}$ and
$\nu(dy) = C_{\alpha} |y|^{-1-\alpha} dy$ is the L\'evy measure. Note that
the compensation term now correctly contains the factor $h(p)$. The
stationary distribution no longer has a simple closed form. Its qualitative
features, in particular the dispersion of the reflected scaled target,
are studied numerically in Section~\ref{sec:numerical}. We do not attempt
to characterize analytically the boundary exponents of the stationary
density for $\alpha<2$; the interaction between the nonlocal jump operator
and the state-dependent diffusion coefficient $h(p)$ makes such a
characterization delicate, and it is not needed for the numerical analysis
of the work.

\subsection{Deterministic Skeleton and Multiple Steady States}\label{sec:escape}

The theory of random perturbations of dynamical systems provides a
classical framework for studying stability, attraction domains, and exit
problems \cite{freidlin2012}. In the present model, the natural object is
the deterministic skeleton obtained by setting $\sigma=0$ and imposing the
quasi-steady-state condition $s=p=\eta g(k)$. Under this condition the
steady states satisfy
\begin{equation}\label{eq:steady_state}
	\eta g(k) k^{\theta - 1} = 1.
\end{equation}
Depending on the parameters, equation \eqref{eq:steady_state} can have one,
two, or three positive solutions. In the multiple-steady-state regime, the
low-capital steady state $k_L$ and the high-capital steady state $k_H$ are
locally stable for the one-dimensional reduced equation
$\dot k = \mu[\eta g(k) k^\theta - k]$, while the intermediate steady state
$k_M$ is unstable. The basin of attraction of $k_H$ is $(k_M, \infty)$, and
that of $k_L$ is $(0, k_M)$. A complete stability analysis of the full
three-dimensional system \eqref{scaledk}--\eqref{scaledp} without
the quasi-steady-state assumption would require linearizing the
three-dimensional vector field at the equilibria $(k_i, \eta g(k_i),
\eta g(k_i))$, $i\in\{L,M,H\}$,  can be referenced from Yuan et al.\
\cite{yuan2026}.

In the stochastic model \eqref{scaledk}--\eqref{scaledp}, the economy can transition from the basin of
$k_H$ to that of $k_L$ if a sufficiently large negative shock drives $k$
below $k_M$. Conversely, a sufficiently large positive shock can drive $k$
from the low-saving basin into the high-saving basin. The frequency of such
transitions depends on the tail index $\alpha$: heavier tails (smaller
$\alpha$) imply more frequent large jumps, and hence plausibly higher
finite-horizon escape frequencies.

In the full coupled system \eqref{scaledk}--\eqref{scaledp}, capital is not a one-dimensional diffusion.
Randomness enters through the scaled saving target $p$, is transmitted to
the scaled saving control $s$, and only then affects $k$. Consequently,
the classical one-dimensional scale-density formula for first-exit
probabilities does not apply directly. For the L\'evy-driven system \eqref{scaledk}--\eqref{scaledp},
first-passage probabilities are determined by the integro-differential
generator or computed numerically. In Section~\ref{sec:numerical} we
therefore report finite-horizon Monte Carlo first-passage probabilities
\begin{equation*}
	q_T(k_0;\alpha) = \Prob_{k_0}\big(\tau_H<\tau_L,\ \tau_H\le T\big),
\end{equation*}
where $\tau_H$ and $\tau_L$ are the first hitting times of prescribed
high- and low-capital targets, rather than closed-form escape
probabilities. If a one-dimensional diffusion approximation is desired, it
must be derived explicitly as a reduced (adiabatic) model and stated as an
approximation.

\subsection{The Role of the Tail Index}\label{sec:sensitivity}

The tail index $\alpha$ plays a dual role in the model \eqref{scaledk}--\eqref{scaledp}. On the one hand, it
governs the frequency and magnitude of jumps in the scaled saving target
$p$. On the other hand, it affects the dispersion of the realized scaled
saving rate $s$ and the finite-horizon transition statistics between
capital regimes.

To isolate the effect of $\alpha$, consider the limiting cases. As
$\alpha\to2$, the L\'evy noise converges to Brownian motion (with the
factor $\sqrt{2}$ as in Section~\ref{sec:stationary}), and the model reduces
to a stochastic Solow model with state-dependent saving and Gaussian
shocks. In this limit, the sample paths of $p$ are continuous, and large
deviations from the mean occur only through the accumulation of many small
shocks. As $\alpha\to1$, the noise becomes extremely heavy-tailed, and the
process $p$ exhibits frequent large jumps. For $\alpha\le 1$, the mean of
$L^{\alpha}(t)$ is not defined, and the model's behavior is qualitatively
different: the scaled saving target can jump arbitrarily far from its mean,
and any notion of a stationary distribution must rely on the reflector
$\mathcal R_{[0,\eta]}$ to keep the process bounded.

For $\alpha\in(1,2)$, the model interpolates between these extremes. The
influence of $\alpha$ on the stationary dispersion of the scaled saving
target and on the realized scaled saving rate is examined numerically in
Section~\ref{sec:numerical}. In the fixed-capital experiments reported
there, the empirical stationary variance of the reflected scaled saving
target decreases as $\alpha$ increases over the tested range. This is a
numerical finding for the specified reflected discrete scheme and
parameter values, not a general monotonicity theorem for all reflected
L\'evy-driven systems.

\subsection{Three-Dimensional Stochastic Dynamics}\label{sec:3d}

Let
$$
X_t=(k_t,s_t,p_t)^\top,
\qquad
h(p)=\sigma\sqrt{p(\eta-p)}.
$$
For the \emph{unreflected} continuous-time system
\eqref{scaledk}--\eqref{scaledp}, the formal infinitesimal
generator is
\begin{equation}
	\begin{aligned}
		\mathcal L_\alpha\phi
		={}&
		\mu(sk^\theta-k)\partial_k\phi+\beta(p-s)\partial_s\phi
		+a(\eta g(k)-p)\partial_p\phi \\
		&+\int_{\R\setminus\{0\}}\Big[
		\phi(k,s,p+h(p)z)-\phi(k,s,p)
		-h(p)z\partial_p\phi\,\mathbf1_{\{|z|<1\}}\Big]\nu_\alpha(dz),
	\end{aligned}
	\label{eq:gen3d}
\end{equation}
where
$$
\nu_\alpha(dz)=C_\alpha|z|^{-1-\alpha}\,dz,
\qquad 0<\alpha<2.
$$
The first three terms describe capital accumulation, saving adjustment,
and mean reversion of the scaled saving target. The integral term
represents the nonlocal contribution of L\'evy jumps. The numerical
scheme applies mirror reflection to the scaled saving target after each
time step. Accordingly, the generator \eqref{eq:gen3d} should be
distinguished from that of a rigorously defined reflected L\'evy process;
the discrete reflected scheme \eqref{eq:disc_p}--\eqref{eq:disc_k} is the
operative model. The full state vector $(k,s,p)$ is used in the Monte
Carlo analysis of joint distributions, capital occupation, and
first-passage events in Section~\ref{sec:numerical}.

\section{Numerical Simulations and Economic Implications}\label{sec:numerical}

This section presents numerical simulations of the discrete-time system
\eqref{eq:disc_p}--\eqref{eq:disc_k} to illustrate the theoretical
results of Section~\ref{sec:stochastic} and to explore the model's
implications for the relationship between state-dependent saving and
macroeconomic volatility.

\subsection{Simulation Methodology}\label{sec:methodology}

The discrete-time system \eqref{eq:disc_p}--\eqref{eq:disc_k} is simulated
over a total horizon of $T=10{,}000$ quarters, corresponding to 2,500
years. With time step $\Delta t=0.05$, one quarter contains $20$ numerical
updates and each trajectory contains $200{,}000$ numerical updates. To
generate the $\alpha$-stable noise terms $\varepsilon_n$, we use the
Chambers--Mallows--Stuck (CMS) algorithm \cite{chambers1976,weron1996},
which produces exact simulations of symmetric $\alpha$-stable random
variables. The algorithm is based on the following representation: if $V$
is uniform on $(-\pi/2,\pi/2)$ and $W$ is exponential with mean 1, then
\begin{equation*}
	\varepsilon = \frac{\sin(\alpha V)}{(\cos V)^{1/\alpha}}
	\left( \frac{\cos((1-\alpha)V)}{W} \right)^{(1-\alpha)/\alpha}
\end{equation*}
is a standard symmetric $\alpha$-stable random variable for $\alpha \neq 1$.
For $\alpha=1$, the symmetric representation reduces to
$$
Z=\tan V,
\qquad
V\sim U\left(-\frac{\pi}{2},\frac{\pi}{2}\right).
$$

The baseline parameters are selected to satisfy the model restrictions and
to illustrate the stochastic dynamics. The baseline parameter values are reported in
Table~\ref{tab:params}.

\begin{table}[htbp]
	\centering
	\caption{Baseline parameter values.}
	\label{tab:params}
	\begin{tabular}{lll}
		\toprule
		Parameter & Description & Value \\
		\midrule
		$\theta$ & Capital share & 0.33 \\
		$s_1$ & Low saving bound (logistic) & 0.05 \\
		$s_2$ & High saving bound (logistic) & 0.25 \\
		$\mu$ & Depreciation rate per unit time & 0.02 \\
		$\gamma$ & Logistic sensitivity & 5 \\
		$\beta$ & Saving adjustment speed & 2 \\
		$a$ & Target adjustment speed & 1 \\
		$\sigma$ & Noise scale & 0.1 \\
		$\alpha$ & Tail index & 1.5 \\
		\bottomrule
	\end{tabular}
\end{table}

The capital share $\theta = 0.33$ is standard in the growth literature. The
saving bounds $s_1 = 0.05$ and $s_2 = 0.25$ are consistent with empirical
estimates of saving rates in low- and high-income countries. The
depreciation rate $\mu = 0.02$ corresponds to an annual rate of
approximately $8\%$ (when $\Delta t$ is measured in years; here it is the
per-unit-time rate in the continuous model, which we discretize with
$\Delta t=0.05$). The logistic sensitivity $\gamma = 5$ implies that the
saving rate transitions from $s_1$ to $s_2$ over a capital range of
approximately $4/\gamma = 0.8$ units around the reference level $k = 1$.
The adjustment speeds $\beta = 2$ and $a = 1$ satisfy
$\beta\Delta t = 0.1$ and $a\Delta t = 0.05$, so the saving control and
target evolve on a faster time scale than capital, as assumed in the
time-scale separation. The noise scale $\sigma = 0.1$ and tail index
$\alpha = 1.5$ are chosen to produce moderate but non-negligible stochastic
fluctuations.

The numerical experiments are based on the scheme
\eqref{eq:disc_p}--\eqref{eq:disc_k}. Time is measured in quarters, with
numerical step $\Delta t=0.05$; thus each quarter contains $20$ numerical
updates. The total simulation horizon is $T=10{,}000$ quarters,
corresponding to $200{,}000$ numerical steps per trajectory. For the
long-run distribution and occupation analyses, the first $5{,}000$ quarters
are discarded as burn-in. Post-burn-in observations are recorded every
four quarters, giving $1{,}250$ sampled observations per trajectory.
Unless otherwise specified, the initial capital is sampled from
$$
k_0\sim U(K_L,K_H),
\qquad
s_0=p_0=\eta g(k_0).
$$
Symmetric $\alpha$-stable innovations are generated using the
Chambers--Mallows--Stuck algorithm. After each scaled saving-target update,
the mirror-reflection map $\mathcal R_{[0,\eta]}$ ensures
$p_n\in[0,\eta]$. Different experiments use different Monte Carlo sample
sizes, which are specified in the corresponding subsections.

\subsection{Stationary Distributions}\label{sec:stationary_dist}

\begin{figure}[H]
	\centering
	\includegraphics[width=1\textwidth]{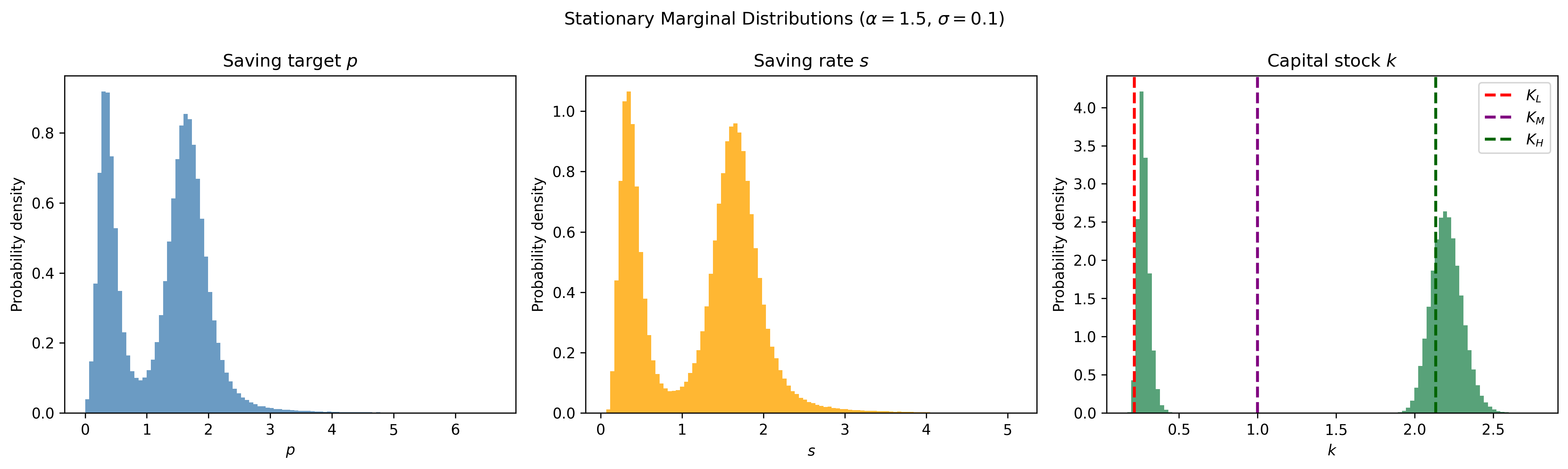}
	\caption{Simulated stationary distributions for the baseline parameters
		($\alpha=1.5$, $\sigma=0.1$). The three panels display the empirical
		probability densities of the scaled saving target $p$, the scaled saving
		control $s$, and the capital stock $k$ after the burn-in period. In the
		right panel, the dashed lines denote the deterministic equilibria $K_L$,
		$K_M$, and $K_H$.}
	\label{fig:stationary}
\end{figure}

Figure~\ref{fig:stationary} presents the empirical long-run distributions
of the scaled saving target $p$, the scaled saving control $s$, and the
capital stock $k$ under the baseline parameters. The capital distribution
exhibits two distinct concentration regions near the deterministic
equilibria $K_L$ and $K_H$, while relatively little probability mass is
observed around the intermediate equilibrium $K_M$. The distributions of
$p$ and $s$ also display two main peaks, reflecting the nonlinear coupling
between the scaled saving target, saving adjustment, and capital
accumulation. These results provide numerical evidence of two dominant
capital-related concentration regions in the simulated system \eqref{eq:disc_p}--\eqref{eq:disc_k}.

The distribution of $s$ is similar to that of $p$, reflecting the
adjustment of the scaled saving control toward the stochastic scaled
saving target. The distribution of $k$ is skewed to the right, with a long
upper tail reflecting the possibility of sustained high-saving episodes.

\begin{figure}[H]
	\centering
	\includegraphics[width=0.77\textwidth]{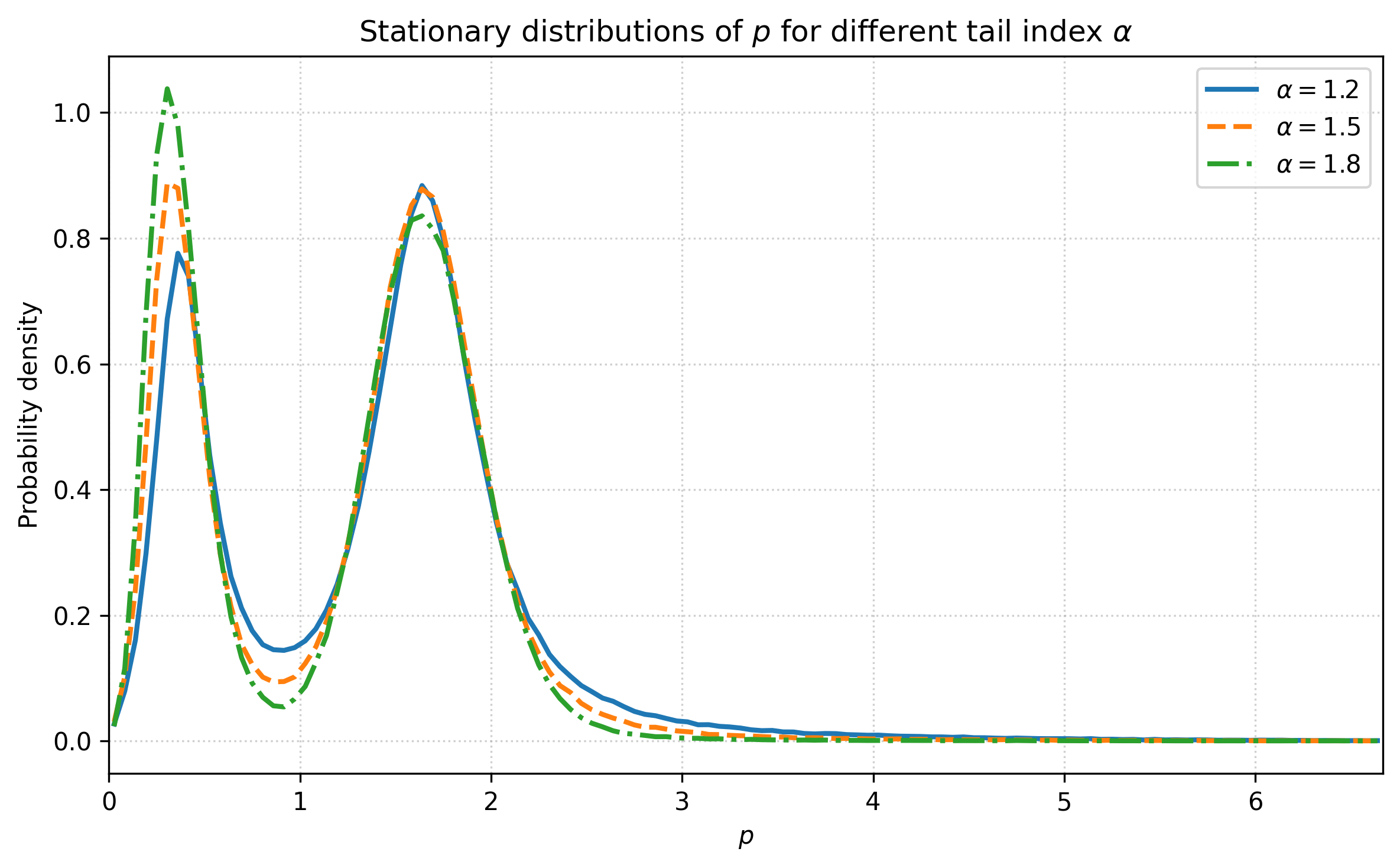}
	\caption{Stationary distributions of the scaled saving target $p$ for
		$\alpha = 1.2, 1.5, 1.8$ in the full coupled system.}
	\label{fig:alpha}
\end{figure}

To illustrate the effect of the tail index $\alpha$,
Figure~\ref{fig:alpha} compares the stationary distributions of $p$ for
$\alpha = 1.2$, $\alpha = 1.5$, and $\alpha = 1.8$. For the tested values
of $\alpha$, smaller tail indices are associated with greater dispersion
of the simulated scaled saving target.

Although $k$ and $s$ evolve endogenously in the full coupled system,
Figure~\ref{fig:alpha} reports only the marginal stationary distribution
of $p$. The effects of $k$ and $s$ are therefore incorporated implicitly
through the joint dynamics rather than displayed separately. In the
full-system simulation, the initial capital stock is sampled as
$k_0 \sim U(K_L,K_H)$, which provides a neutral initialization across the
two main capital regimes.  After initialization, $k_n$ evolves
according to the endogenous capital dynamics, so its subsequent values are
determined by the coupled system rather than by the initial uniform
distribution.

To isolate the direct effect of the L\'evy tail index $\alpha$ on the
scaled saving target $p$, we fix the capital stock and remove the feedback
from $k$ to $p$. Although any deterministic equilibrium could in principle
be used as a reference point, $K_L$ and $K_H$ correspond to the low and
high capital regimes and would therefore introduce regime-specific effects.
We instead choose the intermediate equilibrium $K_M$, which provides a
common reference state between the two regimes. Under the baseline
parametrization, $K_M=1$, and hence
$$
g(K_M)=\frac{s_1+s_2}{2}, \qquad \eta g(K_M)=1.
$$
The fixed-$k$ experiment therefore studies the response of $p$ to
different values of $\alpha$ around the common reference level $p^{\ast}=1$.

Table~\ref{tab:fixed_k_alpha_stats_p} reports the empirical long-run
moments of the scaled saving target for different values of the L\'evy
stability index, with capital fixed at $k=K_M$. The estimated variance
decreases from $0.290433$ at $\alpha=1.2$ to $0.084062$ at $\alpha=1.8$.
This indicates that, under the specified parameters and numerical
reflection mechanism, smaller values of $\alpha$ are associated with
greater dispersion of the scaled saving target. The result is a numerical
observation for the fixed-capital model, rather than a general
monotonicity theorem. The observed decrease in variance is a numerical
result for the fixed-capital experiment under the specified parameter
values and mirror-reflection mechanism; it should not be interpreted as a
general monotonicity property for all reflected L\'evy-driven systems. In
the fixed-capital experiment, we consider the scaled saving target $p_n$
and the scaled saving rate $s_n$ driven by symmetric $\alpha$-stable
L\'evy noise, with $k=K_M$. For $0<\alpha<2$, the driving $\alpha$-stable
L\'evy noise has no finite second moment,
$\E[|L_t^\alpha|^2]=\infty$, but the mirror-reflection mechanism ensures
$0\le p_n\le\eta$. The scaled saving rate is updated by
$s_{n+1}=s_n+\beta\Delta t(p_{n+1}-s_n)$; since $\beta\Delta t=0.1$, this
becomes $s_{n+1}=0.9s_n+0.1p_{n+1}$. Thus both $p_n$ and $s_n$ remain
bounded in $[0,\eta]$ and have finite second moments.

\begin{table}[H]
	\centering
	\caption{Stationary statistics of the scaled saving target $p$ for
		different tail indices $\alpha$ with fixed $k=K_M$.}
	\label{tab:fixed_k_alpha_stats_p}
	\begin{tabular}{cccc}
		\toprule
		$\alpha$ & Mean$(p)$ & Var$(p)$ & Std$(p)$ \\
		\midrule
		1.2 & 1.119489 & 0.290433 & 0.538918 \\
		1.5 & 1.046905 & 0.152357 & 0.390330 \\
		1.8 & 1.012568 & 0.084062 & 0.289935 \\
		\bottomrule
	\end{tabular}
\end{table}

Using post-burn-in Monte Carlo samples, we estimate the stationary
variances by
$$
\widehat{\operatorname{Var}}(p)
=\frac{1}{N}\sum_{i=1}^{N}(p_i-\bar p)^2,
\qquad
\widehat{\operatorname{Var}}(s)
=\frac{1}{N}\sum_{i=1}^{N}(s_i-\bar s)^2.
$$
We further define the variance ratio
$R(\alpha)=\widehat{\operatorname{Var}}(s)/\widehat{\operatorname{Var}}(p)$,
which measures the attenuation of fluctuations from the scaled saving
target $p$ to the scaled saving rate $s$. These quantities represent
empirical variances of the bounded reflected state variables, rather than
the theoretical variance of the driving L\'evy noise. The sample mean in
Table~\ref{tab:fixed_k_alpha_stats_p} exceeds the reference level $p^{\ast}=1$;
this upward bias reflects the asymmetric placement of the reflector
boundaries $0$ and $\eta=2/(s_1+s_2)$ relative to $p^{\ast}=1$ in the
baseline parametrization, and it is a feature of the reflected numerical
scheme rather than of the underlying continuous-time model.

\begin{figure}[H]
	\centering
	\includegraphics[width=0.77\textwidth]{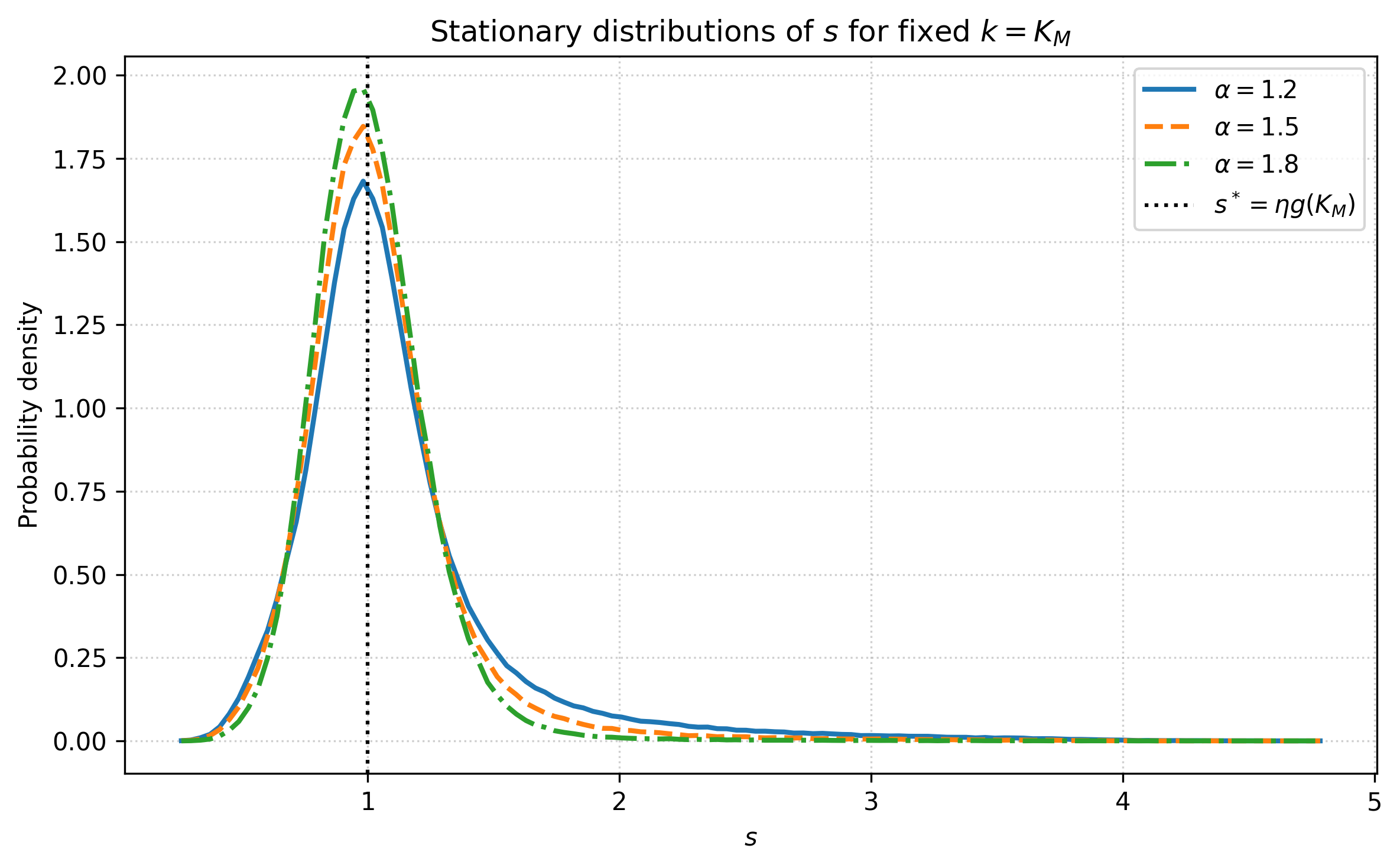}
	\caption{Stationary distributions of the scaled saving control $s$ for
		$\alpha = 1.2, 1.5, 1.8$ with capital fixed at $k=K_M$. The reference
		level $K_M$ is obtained numerically as the intermediate root of the
		deterministic equilibrium condition $k^{1-\theta}=\eta g(k)$, which yields
		$K_M\approx 1.0$. Consequently, the deterministic equilibrium for the
		scaled saving rate is $s^{\ast}=\eta g(K_M)=1$ (indicated by the dotted
		line).}
	\label{fig:alpha_s}
\end{figure}

Figure~\ref{fig:alpha_s} shows the stationary distributions of the scaled
saving control $s$ for $\alpha=1.2,1.5,1.8$, with capital fixed at
$k=K_M$. As $\alpha$ increases from $1.2$ to $1.8$, the distribution
becomes more concentrated around the equilibrium $s^{\ast}=\eta g(K_M)$,
with a higher peak and a lighter right tail. In contrast, smaller values
of $\alpha$ produce a more dispersed distribution with a heavier right
tail. Table \ref{tab:fixed_k_alpha_stats_s} presents the stationary statistics of the scaled saving control $s$ for
		different tail indices $\alpha$ with fixed $k=K_M$. Table \ref{tab:variance_comparison} shows the comparison of the stationary variances of $p$ and $s$ for
		different tail indices $\alpha$. This suggests that the influence of $\alpha$ on $p$ is transmitted
to $s$, consistent with the mechanism $\alpha\to p\to s$.

\begin{table}[H]
	\centering
	\caption{Stationary statistics of the scaled saving control $s$ for
		different tail indices $\alpha$ with fixed $k=K_M$.}
	\label{tab:fixed_k_alpha_stats_s}
	\begin{tabular}{cccc}
		\toprule
		$\alpha$ & Mean$(s)$ & Var$(s)$ & Std$(s)$ \\
		\midrule
		1.2 & 1.118843 & 0.190874 & 0.436891 \\
		1.5 & 1.047146 & 0.101774 & 0.319021 \\
		1.8 & 1.012270 & 0.056470 & 0.237633 \\
		\bottomrule
	\end{tabular}
\end{table}

\begin{table}[H]
	\centering
	\caption{Comparison of the stationary variances of $p$ and $s$ for
		different tail indices $\alpha$.}
	\label{tab:variance_comparison}
	\begin{tabular}{c|cc}
		\toprule
		$\alpha$ & $\mathrm{Var}(p)$ & $\mathrm{Var}(s)$ \\
		\midrule
		1.2 & 0.290433 & 0.190874 \\
		1.5 & 0.152357 & 0.101774 \\
		1.8 & 0.084062 & 0.056470 \\
		\bottomrule
	\end{tabular}
\end{table}

\begin{figure}[H]
	\centering
	\includegraphics[width=0.77\textwidth]{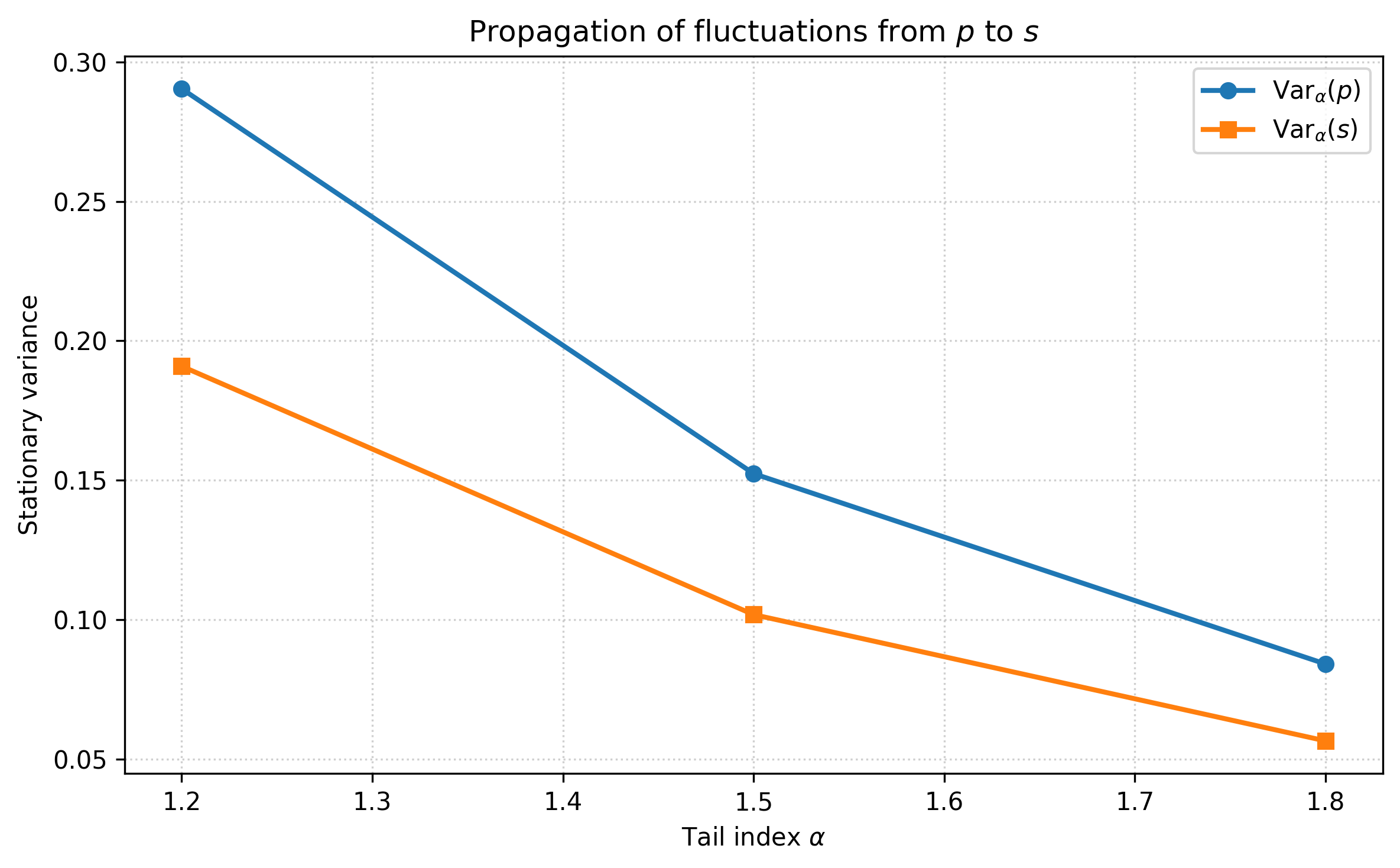}
	\caption{Propagation of fluctuations from the scaled saving target $p$ to
		the scaled saving control $s$.}
	\label{fig:p_s}
\end{figure}

The scaled saving rate $s$ adjusts gradually towards the scaled saving
target $p$ according to $s_{n+1}=s_n+\beta\Delta t(p_{n+1}-s_n)$. In the
baseline setting $\beta=2$ and $\Delta t=0.05$, so
$s_{n+1}=0.9s_n+0.1p_{n+1}$. Hence the dynamics of $s$ act as a smoothing
mechanism: shocks to $p$ are transmitted to $s$, but only partially at
each time step, resulting in attenuated fluctuations in $s$ relative to
$p$, as illustrated in Figure \ref{fig:p_s}.

\begin{figure}[H]
	\centering
	\includegraphics[width=0.77\textwidth]{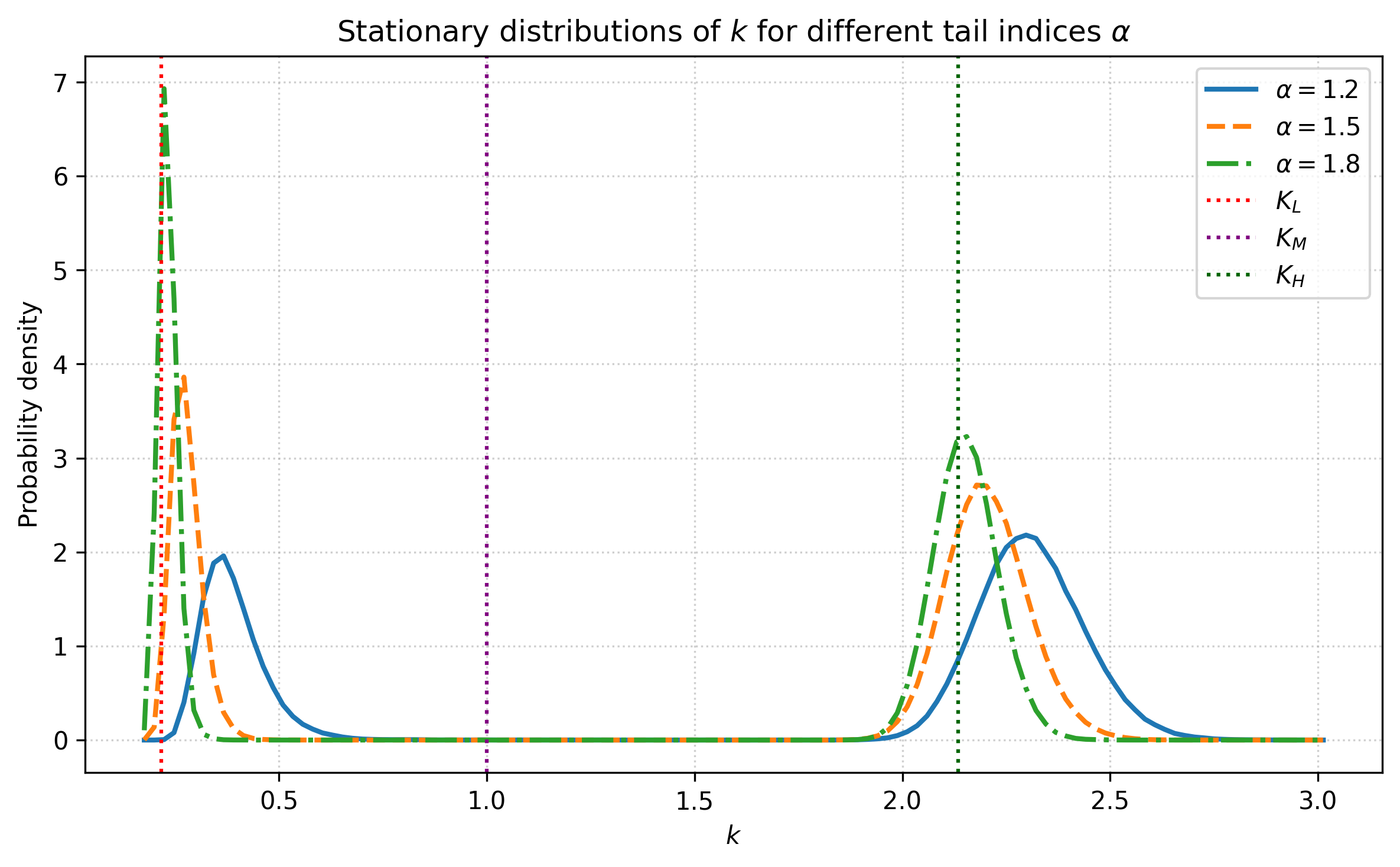}
	\caption{Stationary distributions of $k$ for $\alpha = 1.2, 1.5, 1.8$ in
		the full coupled system. The initial capital stock is drawn from a uniform
		distribution $k_0\sim U(K_L,K_H)$, with initial saving values set to
		$p_0=s_0=\eta g(k_0)$. Vertical dotted lines denote the deterministic
		equilibria $K_L$, $K_M$, and $K_H$.}
	\label{fig:alpha_k}
\end{figure}

Using the discrete-time approximation
\eqref{eq:disc_k}, with $\mu=0.02$ and $\Delta t=0.05$,
$$
k_{n+1}=k_n+0.001\,(s_{n+1}k_n^{\theta}-k_n).
$$
Figure \ref{fig:alpha_k} plots  the stationary distributions of $k$ for $\alpha = 1.2, 1.5, 1.8$ in
		the full coupled system \eqref{scaledk}--\eqref{scaledp}. This shows that changes in the scaled saving rate $s$ are transmitted to
capital $k$, but the adjustment of $k$ is much slower and occurs through
gradual accumulation rather than immediate one-to-one changes.

\subsection{Monte Carlo Joint and Marginal Densities}\label{sec:MC}

Monte Carlo methods provide a numerical framework for studying stochastic
systems driven by L\'evy processes. Related numerical approaches, including
multilevel Monte Carlo methods based on jump-adapted Euler schemes, have
been investigated by Dereich and Li \cite{dereich2016}. For the Monte Carlo
simulation we use the scheme \eqref{eq:disc_p}--\eqref{eq:disc_k} with
$\Delta t=0.05$ and $\alpha_d=a\Delta t$, $\beta_d=\beta\Delta t$,
$\mu_d=\mu\Delta t$.

To investigate the long-run joint behavior of the coupled stochastic
system, we estimate the time-averaged joint probability density of the
capital stock $k$, the scaled saving rate $s$, and the scaled saving
target $p$ using Monte Carlo simulation. The simulation consists of
$N_{\mathrm{MC}}=500$ sample paths over $10{,}000$ quarters, with time step
$\Delta t=0.05$. The first $5{,}000$ quarters are discarded as burn-in, and
the remaining observations are sampled every four quarters, yielding
$1{,}250$ observations per path and $625{,}000$ joint observations in
total. The time-averaged joint density is estimated from post-burn-in
Monte Carlo observations using a three-dimensional histogram. The state
space is discretized into $41\times41\times61$ bins along the $k$, $s$,
and $p$ directions, respectively. For each grid cell the empirical joint
density is
\begin{equation*}
	\widehat{\rho}_{ijk}
	=\frac{N_{ijk}}{N_{\mathrm{obs}}\Delta k\Delta s\Delta p},
\end{equation*}
where $N_{ijk}$ denotes the number of observations in the corresponding
cell and $N_{\mathrm{obs}}=625{,}000$. Only cells with estimated density
above $10^{-4}$ times the maximum density are displayed in the
three-dimensional figure. This threshold is used solely for visualization
and does not modify the underlying density estimates.

\begin{figure}[H]
	\centering
	\includegraphics[width=0.77\textwidth]{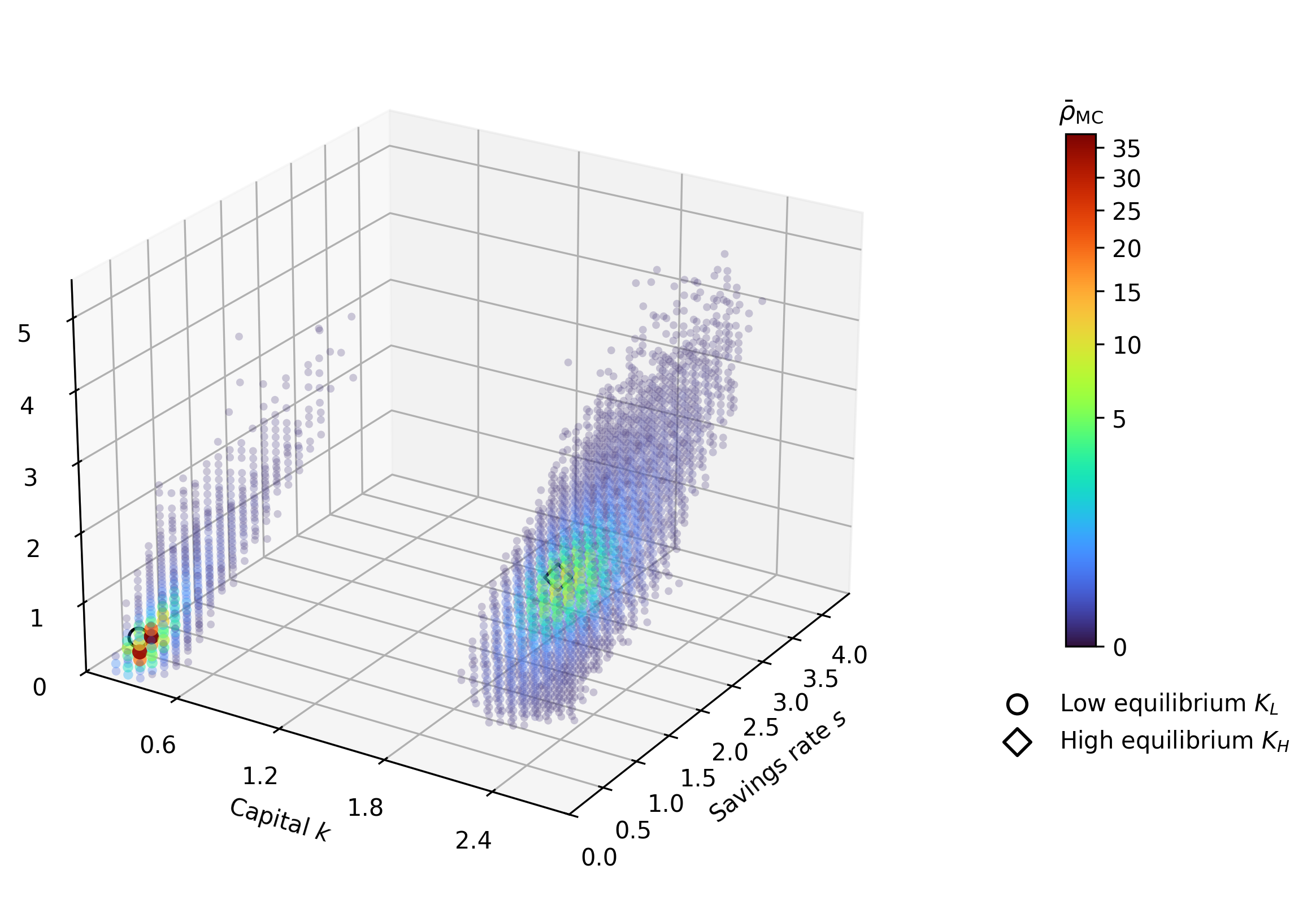}
	\caption{Time-averaged three-dimensional joint probability density
		$\bar{\rho}_{\mathrm{MC}}(k,s,p)$ estimated from Monte Carlo simulations
		under the baseline parameters $\alpha=1.5$, $\sigma=0.1$, and $\gamma=5.0$.
		The markers denote the low- and high-capital equilibria $K_L$ and $K_H$,
		respectively. The color scale represents the estimated joint density. Only
		grid cells satisfying $\rho>0$ and
		$\rho\geq0.0001\rho_{\max}$ are displayed.}
	\label{fig:MC_joint_density}
\end{figure}

Figure~\ref{fig:MC_joint_density} shows the time-averaged joint
probability density of $(k,s,p)$ obtained from Monte Carlo simulation. The
density exhibits two distinct concentration regions, corresponding to the
low-capital and high-capital equilibria, $K_L\approx0.217$ and
$K_H\approx2.135$, respectively. Near $K_L$ the density is concentrated
within a relatively narrow region, whereas the distribution around $K_H$
is more dispersed across the three-dimensional state space. The
intermediate region, particularly near the unstable equilibrium $K_M=1$,
has a substantially lower estimated density. These results indicate that
the stochastic system exhibits two main long-run concentration regions,
associated with the low- and high-capital regimes, while the
intermediate-capital region has relatively low occupancy. The observed
bimodal structure is consistent with the broader finding that nonlinear
stochastic dynamics under L\'evy noise can produce multiple concentration
regions \cite{chechkin2004}.

\begin{figure}[H]
	\centering
	\includegraphics[width=1\textwidth]{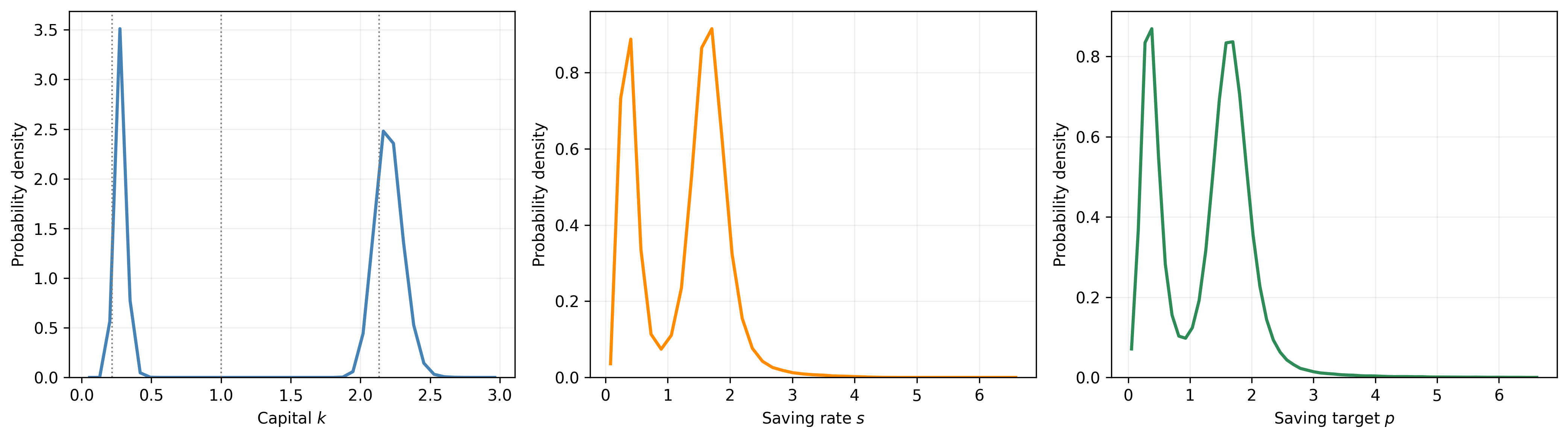}
	\caption{Monte Carlo time-averaged marginal probability densities of
		(a) the capital stock $k$, (b) the scaled saving rate $s$, and (c) the
		scaled saving target $p$, obtained by integrating the three-dimensional
		joint density over the remaining state variables. All three distributions
		exhibit a bimodal structure, consistent with the low- and high-capital
		regimes.}
	\label{fig:MC_marginal}
\end{figure}

To examine the individual state variables, we compute the marginal
densities
\begin{equation*}
	\bar{\rho}_k(k)=\iint\bar{\rho}(k,s,p)\,ds\,dp,
	\quad
	\bar{\rho}_s(s)=\iint\bar{\rho}(k,s,p)\,dk\,dp,
	\quad
	\bar{\rho}_p(p)=\iint\bar{\rho}(k,s,p)\,dk\,ds.
\end{equation*}
The marginal densities are obtained by numerical integration of the
estimated three-dimensional joint density. For example, the capital
marginal is approximated by
$\widehat{\rho}_{k,i}=\sum_{j=1}^{N_s}\sum_{\ell=1}^{N_p}
\widehat{\rho}_{ij\ell}\Delta s\Delta p$. The marginal densities of $s$
and $p$ are computed analogously. Thus, Figure~\ref{fig:MC_marginal}
represents the one-dimensional marginal distributions derived from the
same Monte Carlo joint density shown in Figure~\ref{fig:MC_joint_density}.

Figure~\ref{fig:MC_marginal} presents the Monte Carlo time-averaged
marginal densities of the three state variables $k$, $s$, and $p$.
Figure~\ref{fig:MC_marginal}(a) shows a bimodal distribution of capital
$k$, with two distinct peaks near the stable equilibria $K_L\approx0.217$
and $K_H\approx2.135$. The density near $K_L$ is more concentrated, while
the distribution around $K_H$ is relatively broader. The density is close
to zero in the intermediate region, particularly near the unstable
equilibrium $K_M=1$. Figure~\ref{fig:MC_marginal}(b) shows a similar
bimodal structure in the scaled saving rate $s$, with two peaks at low
and high saving levels. This is consistent with the adjustment of the
scaled saving rate towards the stochastic scaled saving target.
Figure~\ref{fig:MC_marginal}(c) illustrates the marginal density of the
scaled saving target $p$. Two distinct peaks are also observed, indicating
two main concentration regions under the state-dependent L\'evy
perturbations. Overall, all three marginal densities exhibit bimodal
structure, consistent with the two concentration regions identified in
the three-dimensional Monte Carlo joint density. These findings support
the presence of distinct low- and high-capital regimes, with relatively
low long-run occupancy of intermediate states.

\subsection{Long-Run Capital Regime Occupation}\label{sec:occupation}

The long-run behavior of capital is examined by considering the proportion
of time that the stochastic system spends in the low- and high-capital
regimes. Noise-induced transitions between attracting regions have been
studied in L\'evy-driven stochastic systems \cite{imkeller2006meta}. The
intermediate unstable equilibrium $K_M$ is used as a threshold to
distinguish the two regimes:
$$
D_L=\{k:k<K_M\},
\qquad
D_H=\{k:k>K_M\}.
$$
For a given L\'evy stability index $\alpha$, let $\pi_{\alpha,k}(k)$
denote the stationary marginal density of capital in the full coupled
system, provided that such a density exists. The corresponding stationary
occupation probabilities are
\begin{equation*}
	\Pi_L(\alpha) = \int_0^{K_M}\pi_{\alpha,k}(k)\,dk,
	\qquad
	\Pi_H(\alpha) = \int_{K_M}^{\infty}\pi_{\alpha,k}(k)\,dk,
\end{equation*}
with $\Pi_L(\alpha)+\Pi_H(\alpha)=1$ provided the stationary distribution
assigns zero mass to $K_M$.

In the numerical experiments the occupation probabilities are
approximated using Monte Carlo observations after discarding an initial
burn-in period. The initial capital is sampled from $k_0\sim U(K_L,K_H)$,
with the remaining initial variables specified by
$s_0=p_0=\eta g(k_0)$. For each value of $\alpha$, $M=500$ independent
stochastic trajectories are simulated over $T=10{,}000$ quarters, with
time step $\Delta t=0.05$. The first $5{,}000$ quarters are discarded as
burn-in, and capital is subsequently recorded every four quarters. Let
$k_{t_i}^{(j)}$ denote the capital observation at sampling time $t_i$
along trajectory $j$, where $j=1,\ldots,M$ and $i=1,\ldots,N$. The
empirical occupation fractions are estimated by
\begin{equation*}
	\widehat{\Pi}_L(\alpha)=\frac{1}{MN}\sum_{j=1}^{M}\sum_{i=1}^{N}
	\mathbf{1}_{\{k_{t_i}^{(j)}<K_M\}},
	\qquad
	\widehat{\Pi}_H(\alpha)=\frac{1}{MN}\sum_{j=1}^{M}\sum_{i=1}^{N}
	\mathbf{1}_{\{k_{t_i}^{(j)}>K_M\}}.
\end{equation*}
Here $\mathbf{1}_{\{\cdot\}}$ denotes the indicator function. These
estimates measure the fractions of post-burn-in observations belonging to
the two capital regimes. The numerical results are reported in
Table~\ref{tab:occupation}. The $95\%$ confidence intervals are calculated
using path-level Monte Carlo standard errors \cite{asmussen2007}. To
further examine changes between capital regimes, the final column reports
the number of trajectories that cross $K_M$ at least once during the full
simulation horizon.

\begin{table}[H]
	\centering
	\caption{Empirical capital regime occupation for different L\'evy
		stability indices.}
	\label{tab:occupation}
	\small
	\begin{tabular}{c c c c c}
		\toprule
		$\alpha$ & $\widehat{\Pi}_L$ & $\widehat{\Pi}_H$
		& $95\%$ CI for $\widehat{\Pi}_L$ & Paths crossing $K_M$ \\
		\midrule
		1.2 & 0.3083 & 0.6917 & $[0.2678, 0.3487]$ & 49/500 \\
		1.5 & 0.3760 & 0.6240 & $[0.3335, 0.4185]$ & 22/500 \\
		1.8 & 0.3920 & 0.6080 & $[0.3492, 0.4348]$ & 12/500 \\
		\bottomrule
	\end{tabular}
\end{table}

From Table~\ref{tab:occupation}, for $\alpha=1.2$ the estimated
low-capital occupation fraction is $0.3083$, while the high-capital
occupation fraction is $0.6917$. Approximately $30.83\%$ of the
post-burn-in capital observations belong to the low-capital regime,
compared with $69.17\%$ in the high-capital regime. The $95\%$ confidence
interval for $\widehat{\Pi}_L$ is $[0.2678,0.3487]$. Furthermore, $49$ of
$500$ simulated trajectories cross $K_M$ at least once, corresponding to
$9.8\%$ of all trajectories. For $\alpha=1.5$ the low-capital occupation
fraction increases to $0.3760$, while the high-capital occupation fraction
decreases to $0.6240$. The corresponding $95\%$ confidence interval for
$\widehat{\Pi}_L$ is $[0.3335,0.4185]$. In this case $22$ out of $500$
trajectories cross $K_M$ at least once, representing $4.4\%$ of the
simulated paths. For $\alpha=1.8$ the estimated occupation fractions are
$0.3920$ and $0.6080$ for the low- and high-capital regimes, respectively.
The $95\%$ confidence interval for the low-capital occupation fraction is
$[0.3492,0.4348]$. Only $12$ out of $500$ trajectories cross $K_M$,
corresponding to $2.4\%$ of the simulated paths. Comparing the three
cases, the high-capital occupation fraction decreases from $0.6917$ to
$0.6080$ as $\alpha$ increases from $1.2$ to $1.8$, a reduction of $8.37$
percentage points. Meanwhile, the proportion of trajectories crossing
$K_M$ decreases from $9.8\%$ to $2.4\%$. These results suggest that, under
the baseline parameter settings, heavier-tailed L\'evy noise is associated
with a higher empirical occupation fraction of the high-capital regime and
more frequent threshold-crossing events across simulated trajectories.

However, the occupation fractions and crossing counts describe different
aspects of the capital dynamics. The former measure the proportions of
post-burn-in observations in each regime, whereas the latter count
trajectories crossing $K_M$ at least once over the full simulation
horizon, regardless of direction. Moreover, crossing $K_M$ does not
necessarily imply that capital reaches the opposite stable equilibrium.
Since transitions are relatively rare, the empirical occupation fractions
may still depend on the initial capital distribution. They are therefore
interpreted as finite-horizon occupation estimates rather than conclusive
evidence of a unique stationary distribution.

\subsection{First-Passage Probabilities and Restricted Mean First-Passage
	Times}\label{sec:escape_sim}

The previous subsection distinguishes the low- and high-capital regimes
using the intermediate unstable equilibrium $K_M$. However, crossing
$K_M$ does not necessarily mean that capital has reached a region close to
the opposite stable equilibrium. Transition paths in stochastic dynamical
systems driven by stable L\'evy noise have been investigated using
most-probable-path approaches \cite{huang2019}. First-hitting
probabilities between prescribed target sets have been studied using
committor functions \cite{aristoff2024}. Related exit-probability and
mean-exit-time problems for L\'evy-driven systems have also been
investigated in \cite{mao2024}. Here we instead focus on finite-horizon
first-passage probabilities between prescribed capital targets. To examine
capital movements relative to more clearly defined target regions, we
introduce two thresholds
$$
A_L=K_L+\delta(K_H-K_L),
\qquad
A_H=K_H-\delta(K_H-K_L),
$$
where $\delta\in(0,0.5)$ controls the positions of the two targets. In the
baseline experiment we set $\delta=0.2$, giving $A_L\approx0.6006$ and
$A_H\approx1.7512$. These thresholds satisfy
$K_L<A_L<K_M<A_H<K_H$. They are used to define first-passage events,
allowing us to distinguish a crossing of $K_M$ from the first arrival at a
more clearly specified capital target.

\begin{figure}[H]
	\centering
	\includegraphics[width=0.77\textwidth]{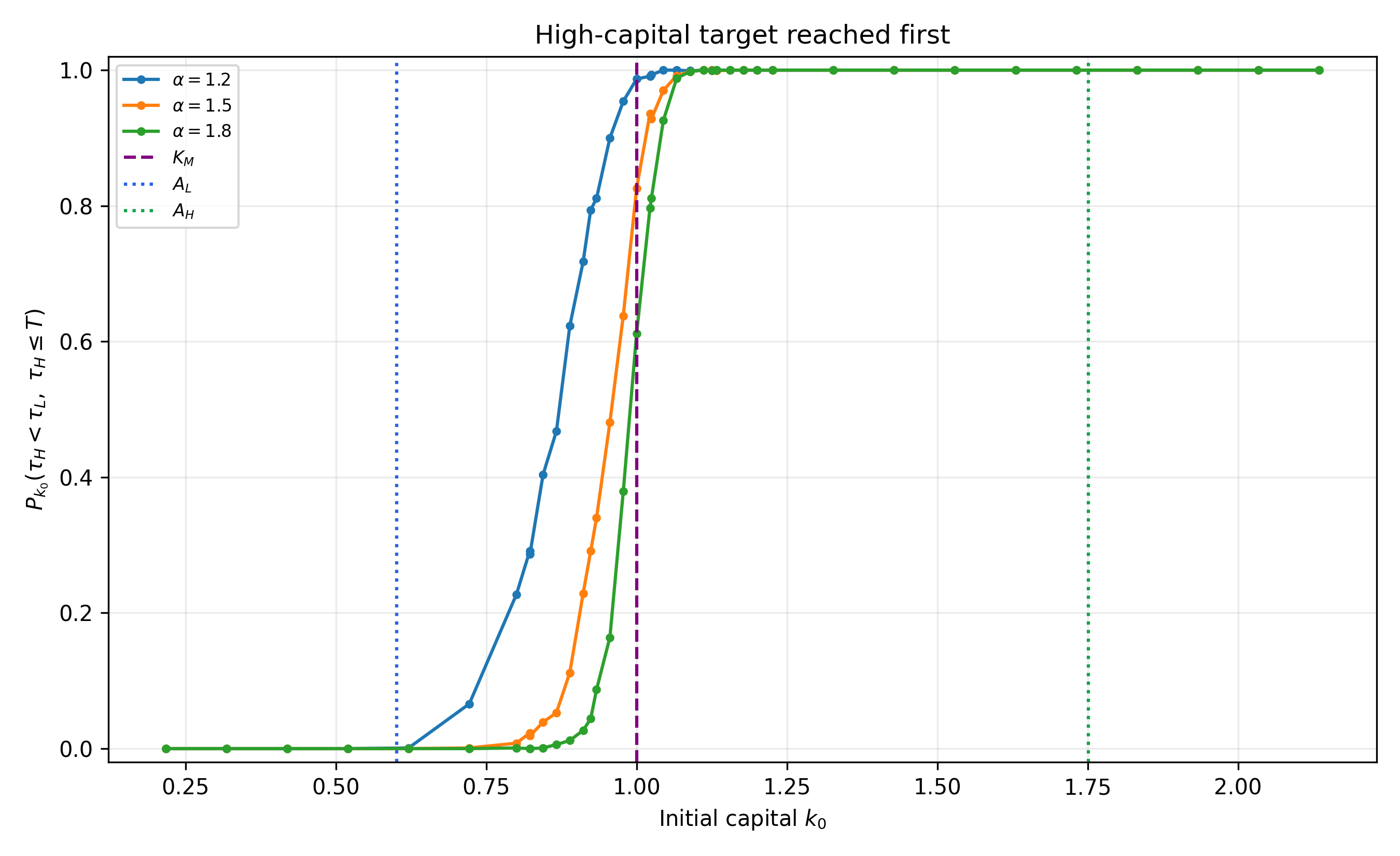}
	\caption{Escape probability
		$q_T(k_0;\alpha)=P_{k_0}(\tau_H<\tau_L,\ \tau_H\leq T)$ for different
		$\alpha$.}
	\label{fig:escape}
\end{figure}

To examine which capital target is reached first, define the first hitting
times of the low- and high-capital targets as
$$
\tau_L=\inf\{t\geq0:k_t\leq A_L\},
\qquad
\tau_H=\inf\{t\geq0:k_t\geq A_H\}.
$$
Comparing these two times allows us to identify which target is reached
first. In particular, the event $\tau_H<\tau_L$ means that capital reaches
the high-capital target before the low-capital target. For a given initial
capital level $k_0$, we define the finite-horizon probability of reaching
the high-capital target first as
$$
q_T(k_0;\alpha)=P_{k_0}(\tau_H<\tau_L,\ \tau_H\leq T),
$$
where $T$ is the observation horizon and the remaining initial variables
are set to $s_0=p_0=\eta g(k_0)$. This probability is estimated using
Monte Carlo simulation. For each $\alpha\in\{1.2,1.5,1.8\}$ a fixed grid
of initial capital values is selected over $[K_L,K_H]$. From each initial
condition, $R=1{,}000$ independent trajectories are simulated. The
estimated probability is
$$
\widehat q_T(k_0;\alpha)=\frac{1}{R}\sum_{r=1}^{R}
\mathbf{1}_{\{\tau_H^{(r)}<\tau_L^{(r)},\,\tau_H^{(r)}\leq T\}}.
$$
Here the indicator equals one if a trajectory reaches the high-capital
target before the low-capital target within the observation horizon, and
zero otherwise. Figure~\ref{fig:escape} plots $\widehat q_T(k_0;\alpha)$
against the initial capital level $k_0$ for different values of $\alpha$.
It shows how the probability of reaching the high-capital target first
varies with both the initial capital level and the L\'evy stability index.
Unlike the occupation fractions in Section~\ref{sec:occupation}, the
first-passage probabilities describe which capital target is reached first
from a given initial condition. Crossing $K_M$ alone does not imply that
capital has reached the high-capital target $A_H$. Furthermore, reaching
$A_H$ does not necessarily indicate a permanent transition to the
high-capital stable equilibrium.

Figure~\ref{fig:escape} shows that the probability of reaching the
high-capital target first increases with the initial capital $k_0$. The
probability curves also depend on the L\'evy tail index $\alpha$. For
$\alpha=1.2$ the curve shifts to the left, indicating a higher probability
of reaching the high-capital target from relatively low initial capital
levels. As $\alpha$ increases, the transition becomes more concentrated
around the intermediate equilibrium $K_M$. These results suggest that
heavier-tailed shocks can influence the selection between the low- and
high-capital regimes. For initial capital levels already within either
target region, the corresponding hitting time is zero; hence the
probabilities at the two ends of the curves are determined by the
prescribed target thresholds.

\begin{figure}[H]
	\centering
	\includegraphics[width=0.9\textwidth]{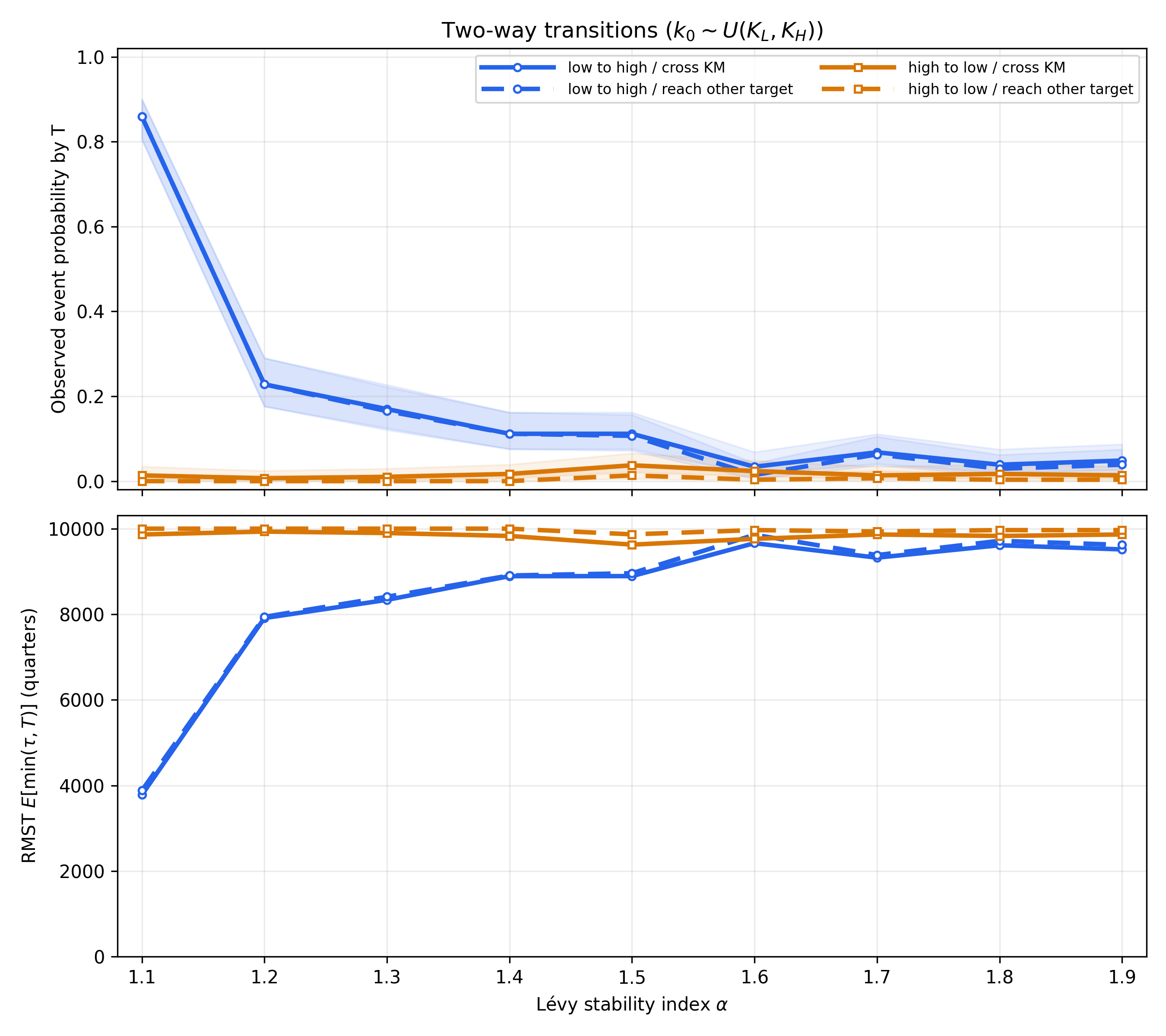}
	\caption{Bidirectional first-passage probabilities and restricted mean
		first-passage times as functions of the L\'evy tail index $\alpha$. RMST
		(Restricted Mean Survival Time) is defined as
		$\mathrm{RMST}_T=\E[\min(\tau,T)]$, where $\tau$ is the first-passage time
		and $T$ is the observation horizon. The shaded bands represent $95\%$
		confidence intervals for the estimated first-passage probabilities
		obtained from Monte Carlo simulations.}
	\label{fig:exit_time}
\end{figure}

Figure~\ref{fig:exit_time} examines bidirectional first-passage events for
different values of $\alpha$. The initial capital is sampled from
$k_0\sim U(K_L,K_H)$, and trajectories are grouped according to whether
$k_0<K_M$ or $k_0>K_M$. We distinguish between the first crossing of $K_M$
and the first arrival at the opposite capital target. The corresponding
event probabilities and restricted mean first-passage times are compared
over the finite observation horizon $T$ ($T=10{,}000$ quarters).

From Figure~\ref{fig:exit_time}, the upper panel shows the probabilities
of crossing the unstable equilibrium $K_M$ and reaching the opposite
capital target within the observation horizon $T$. The blue curves
represent transitions from the low-capital regime to the high-capital
regime, while the orange curves represent transitions in the opposite
direction. For low-to-high transitions, the event probability is highest
at $\alpha=1.1$, reaching approximately $0.86$. As $\alpha$ increases to
$1.2$, the probability decreases sharply to approximately $0.23$. It then
declines more gradually, reaching approximately $0.11$ at $\alpha=1.4$ and
remaining at a similar level at $\alpha=1.5$. At $\alpha=1.6$ the
probability falls to approximately $0.03$, followed by relatively small
fluctuations over $\alpha\in[1.7,1.9]$. These results indicate an overall
decreasing trend in low-to-high first-passage probabilities as $\alpha$
increases, although the relationship is not strictly monotonic. In
particular, the pronounced difference between $\alpha=1.1$ and
$\alpha=1.2$ suggests substantial sensitivity to the tail index in this
parameter range. The solid and dashed blue curves remain close throughout
most of the parameter range. This indicates that the proportions of
trajectories crossing $K_M$ and reaching the high-capital target $A_H$
within $T$ are similar. However, these are distinct first-passage events,
and their similar probabilities do not necessarily imply that the
additional time required to reach $A_H$ after crossing $K_M$ is
negligible. In contrast, the high-to-low probabilities remain close to
zero across the entire range of $\alpha$. The orange curves show only
small fluctuations, with a slight increase in the probability of crossing
$K_M$ around $\alpha=1.5$. The probability of reaching the low-capital
target $A_L$ is generally no greater than the probability of crossing
$K_M$, as reaching $A_L$ requires the trajectory to cross the regime
boundary first. Overall, the upper panel reveals a pronounced asymmetry
between the two transition directions. Under the specified initial
conditions and model parameters, low-to-high first-passage events occur
substantially more frequently than high-to-low events within the
observation horizon.

To summarize first-passage times over a finite observation horizon, we use
the restricted mean first-passage time, which is analogous to the
restricted mean survival time in time-to-event analysis \cite{royston2013}.
The lower panel presents the RMST for the same four events. For a
first-passage time $\tau$, the RMST is defined as
$\operatorname{RMST}(T)=\E[\min(\tau,T)]$. Thus trajectories that reach the
specified target before $T$ contribute their observed first-passage times,
while trajectories that do not reach the target contribute the truncation
time $T$. For low-to-high transitions, the RMST is approximately $3800$
quarters at $\alpha=1.1$. It increases sharply to approximately $7900$
quarters at $\alpha=1.2$ and continues to rise as $\alpha$ increases. By
$\alpha=1.6$ the RMST reaches approximately $9700$ quarters. For larger
values of $\alpha$ the RMST remains close to the observation horizon, with
relatively small fluctuations. This pattern is consistent with the
reduction in low-to-high first-passage probabilities observed in the upper
panel. At smaller values of $\alpha$, a larger proportion of trajectories
reach the opposite regime within $T$, resulting in a lower restricted mean
first-passage time. As $\alpha$ increases, fewer trajectories complete the
event before $T$, causing the RMST to approach the truncation horizon. For
high-to-low transitions, the RMST remains close to $10{,}000$ quarters over
the entire range of $\alpha$. Both orange curves exhibit only limited
variation, suggesting that most trajectories starting in the high-capital
regime do not reach the corresponding opposite-side thresholds within the
observation period. It is important to distinguish the RMST from the
unrestricted mean first-passage time $\E[\tau]$. An RMST close to $T$ does
not imply that trajectories typically complete the transition at time $T$.
Rather, it may reflect a large proportion of trajectories for which no
transition is observed before the simulation terminates. The solid and
dashed curves also distinguish two different stopping times. Since a
trajectory must first cross $K_M$ before reaching the opposite target, we
have $\tau_{\mathrm{cross}}\leq\tau_{\mathrm{target}}$, and consequently
$\operatorname{RMST}_{\mathrm{cross}}(T)\leq
\operatorname{RMST}_{\mathrm{target}}(T)$. The small differences between
the corresponding curves should therefore be interpreted with care,
especially when both RMST values are close to $T$ because of unresolved
trajectories.

The two panels provide complementary information about the stochastic
transitions between capital regimes. The upper panel measures the
probability that a specified first-passage event occurs before $T$,
whereas the lower panel describes the corresponding truncated waiting
time. Let $q_T(\alpha)=P_\alpha(\tau\leq T)$ denote the finite-horizon
first-passage probability. Then
$$
\operatorname{RMST}_{\alpha}(T)
=q_T(\alpha)\E_{\alpha}[\tau\mid\tau\leq T]
+\bigl(1-q_T(\alpha)\bigr)T.
$$
This identity explains the relationship between the two panels. A lower
event probability increases the contribution of trajectories that do not
reach the target before $T$, while the actual arrival times of successful
trajectories also affect the RMST. For example, at $\alpha=1.1$ the
low-to-high crossing probability is approximately $0.86$, accompanied by
an RMST of approximately $3800$ quarters. At $\alpha=1.6$ the
corresponding probability falls to approximately $0.03$, while the RMST
increases to approximately $9700$ quarters. These results are consistent
with the mathematical relationship between first-passage probabilities and
restricted waiting times. A second important observation is the asymmetry
between low-to-high and high-to-low transitions. Although the driving
L\'evy noise is symmetric, the capital dynamics are governed by a
nonlinear coupled system involving capital accumulation, saving
adjustment, and a state-dependent saving target. Consequently, symmetric
stochastic forcing does not necessarily produce symmetric transition
probabilities. From an economic perspective, the numerical results suggest
that the high-capital regime exhibits greater persistence than the
low-capital regime under the considered model specification. In particular,
transitions from the low-capital region toward the high-capital target are
more frequently observed, especially for smaller values of $\alpha$.

These findings indicate that the L\'evy stability index plays an important
role in the finite-horizon transition behavior of the stochastic capital
system. However, the observed asymmetry and parameter sensitivity are
numerical findings for the specified model and initial distributions,
rather than analytical proofs of global stability or monotonicity.

\subsection{Macroeconomic Volatility}\label{sec:volatility}

A central implication of the model is that state-dependent saving can
amplify macroeconomic volatility. To quantify this effect, we compute the
coefficient of variation (CV) of output per effective worker,
$y=k^{\theta}$, for different values of the logistic sensitivity $\gamma$
and the tail index $\alpha$.

\begin{figure}[H]
	\centering
	\includegraphics[width=0.77\textwidth]{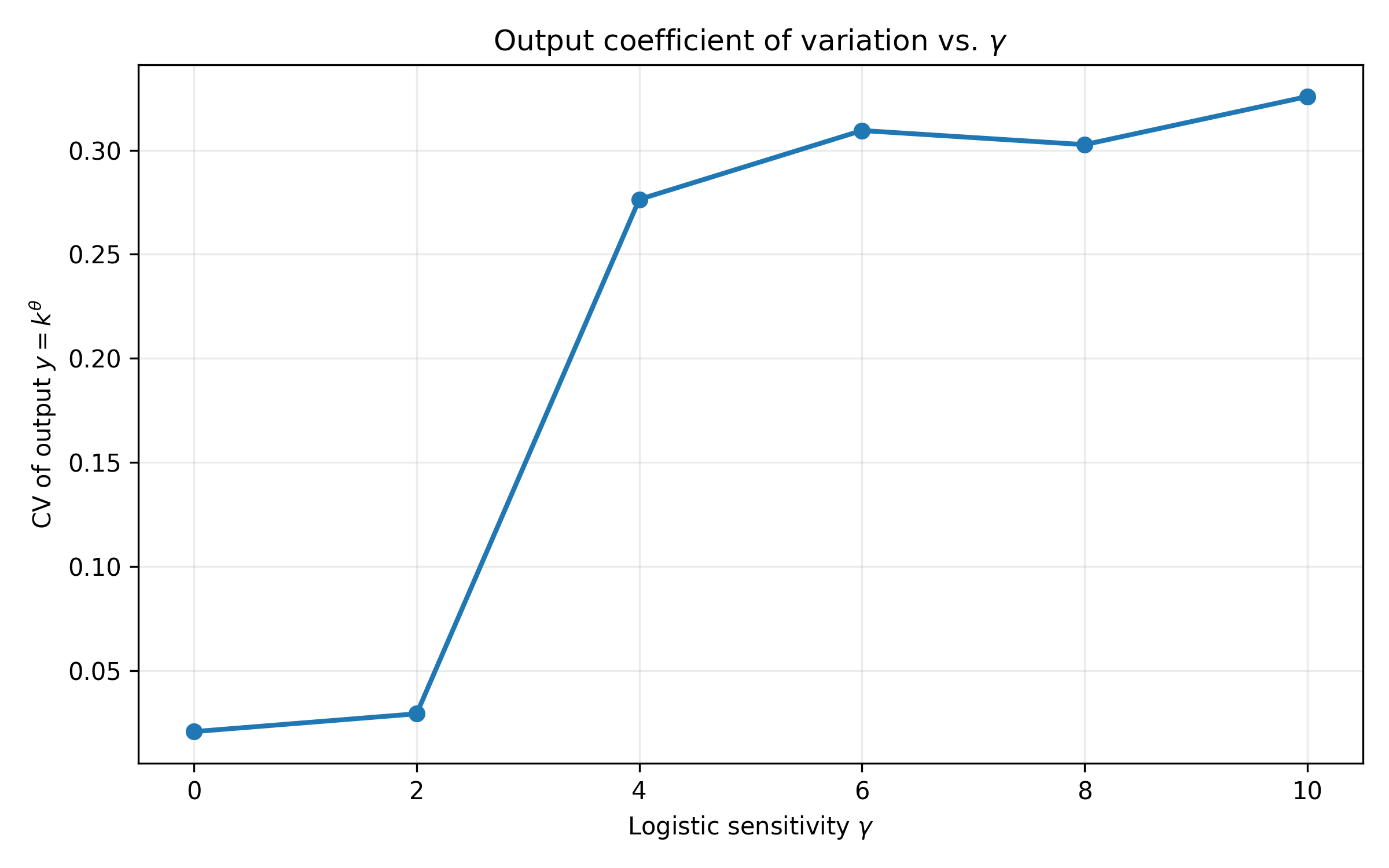}
	\caption{Coefficient of variation of output as a function of $\gamma$ for
		$\alpha = 1.5$.}
	\label{fig:cv_gamma}
\end{figure}

Figure~\ref{fig:cv_gamma} plots the CV of $y$ as a function of $\gamma$
for $\alpha=1.5$. The coefficient of variation is defined by
$$
\mathrm{CV}(y)=\frac{\sqrt{\operatorname{Var}(y)}}{\E[y]},
\qquad
y=k^{\theta},\quad \theta=0.33.
$$
The coefficient of variation is estimated from pooled post-burn-in Monte
Carlo observations. Let $Y_i=k_i^{\theta}$ denote the simulated output
observations. The empirical coefficient of variation is
\begin{equation*}
	\widehat{\mathrm{CV}}(Y)
	=\frac{\sqrt{\frac1N\sum_{i=1}^{N}Y_i^2
			-\left(\frac1N\sum_{i=1}^{N}Y_i\right)^2}}
	{\frac1N\sum_{i=1}^{N}Y_i}.
\end{equation*}
Since the observations are pooled across trajectories and sampling times,
the estimated CV captures both within-regime fluctuations and differences
between the low- and high-capital regimes. As a consequence, a larger
pooled CV may reflect either amplification of fluctuations within regimes
or increased separation between regimes (or both). A decomposition of
total variance into within- and between-regime components would be
required to distinguish these two channels; we regard this as a useful
extension and focus here on the pooled CV as a descriptive summary.

Figure~\ref{fig:cv_gamma} shows that output variability depends
nonlinearly on the saving sensitivity parameter $\gamma$. For $\gamma=0$
the saving target function becomes constant, and the estimated CV is
approximately $0.020$. When $\gamma$ increases to $2$, the CV rises only
slightly to approximately $0.029$, indicating relatively low output
variability in this parameter range. A substantial increase occurs between
$\gamma=2$ and $\gamma=4$, where the CV rises from approximately $0.029$
to $0.277$. This represents the most pronounced change observed in
Figure~\ref{fig:cv_gamma}. The CV remains relatively high, taking values
of approximately $0.309$, $0.302$, and $0.325$ at $\gamma=6$, $8$, and
$10$, respectively. Thus, although output variability generally increases
with $\gamma$, the relationship is not strictly monotonic. The results
suggest that the effect of saving sensitivity on output variability
becomes less pronounced once $\gamma$ exceeds $4$. The nonlinear behavior
observed in Figure~\ref{fig:cv_gamma} may be related to changes in the
deterministic equilibrium structure of the capital accumulation system.
Differentiating the scaled saving target function in
\eqref{eq:SigmoidFunction} with respect to capital gives
$$
g'(k)=\frac{(s_2-s_1)\gamma\exp[-\gamma(k-1)]}
{\left(1+\exp[-\gamma(k-1)]\right)^2}.
$$
At $k=1$ this reduces to
$g'(1)=(s_2-s_1)\gamma/4=0.05\gamma$. Therefore, increasing $\gamma$
strengthens the local sensitivity of the saving target to changes in
capital around $k=1$. The deterministic equilibrium condition is
$k^{1-\theta}=\eta g(k)$. Under the baseline parameter values, $k^{\ast}=1$
is always an equilibrium. The derivative of the equilibrium equation at
this point is
$$
\left.\frac{d}{dk}\left(k^{1-\theta}-\eta g(k)\right)\right|_{k=1}
=(1-\theta)-\frac{\eta(s_2-s_1)\gamma}{4}.
$$
The corresponding critical parameter is
$$
\gamma_c=\frac{4(1-\theta)}{\eta(s_2-s_1)}\approx 2.01.
$$
This value marks a change in the local stability of the central
deterministic equilibrium. Consequently, the sharp increase in the CV
between $\gamma=2$ and $\gamma=4$ occurs across a parameter range
containing this deterministic stability threshold. One possible
explanation is that, for larger $\gamma$, the emergence of distinct stable
capital regimes increases the dispersion of output across simulated
trajectories. However, the link between deterministic stability and output
variability requires further investigation, and further analysis of the
capital distribution is needed to verify this mechanism.

\begin{figure}[H]
	\centering
	\includegraphics[width=0.77\textwidth]{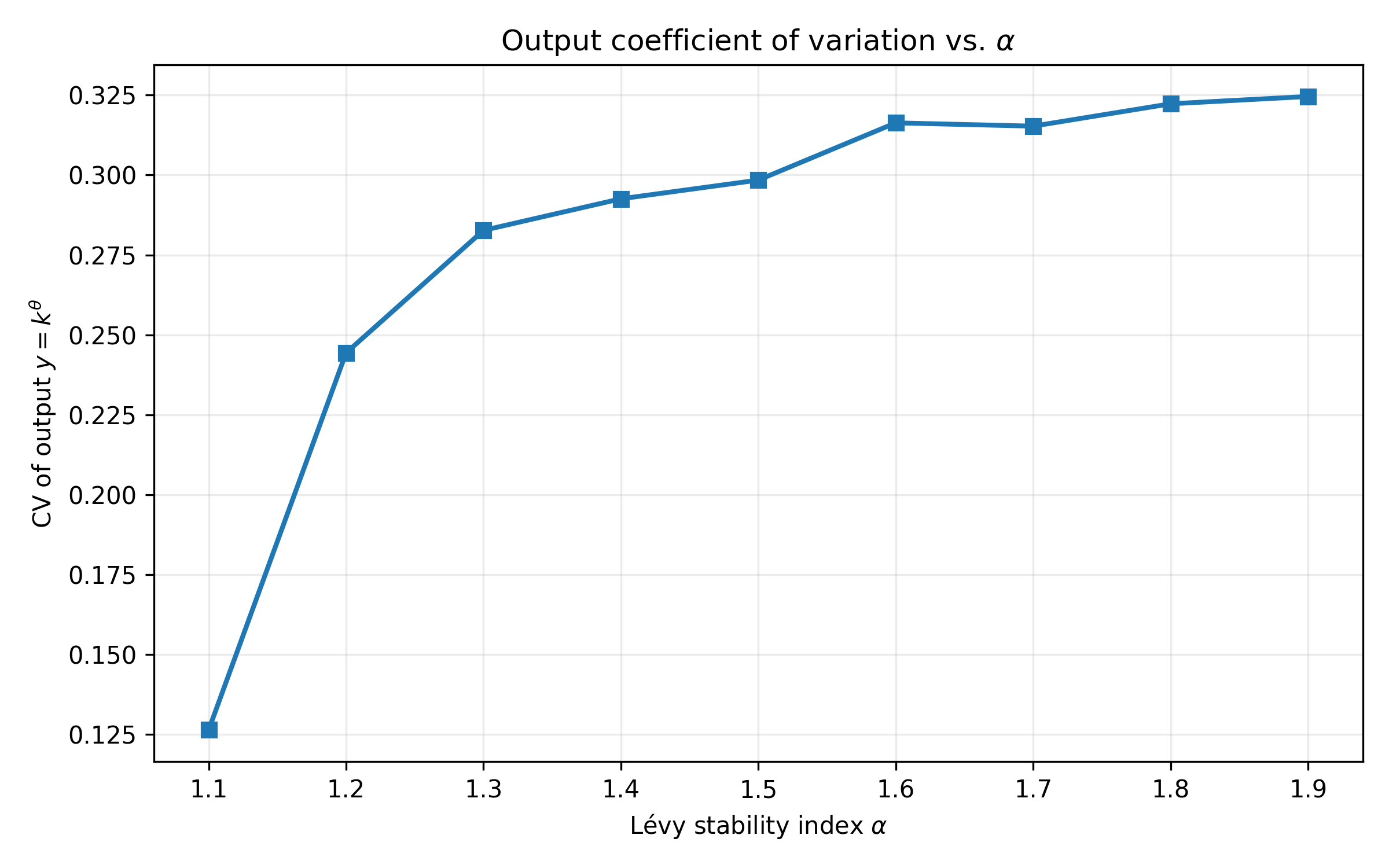}
	\caption{Coefficient of variation of output as a function of $\alpha$ for
		$\gamma = 5$.}
	\label{fig:cv_alpha}
\end{figure}

Figure~\ref{fig:cv_alpha} plots the coefficient of variation of output as
a function of the L\'evy stability index $\alpha$, with $\gamma=5$ and
$\sigma=0.1$ held fixed. The estimated CV increases over the tested
parameter range, indicating that heavier-tailed shocks do not necessarily
generate greater relative output variability.

Figure~\ref{fig:cv_alpha} presents the CV of output as a function of the
L\'evy stability index $\alpha$, with $\gamma=5$ and $\sigma=0.1$ fixed.
The CV generally increases with $\alpha$, rising from approximately
$0.126$ at $\alpha=1.1$ to $0.244$ at $\alpha=1.2$, and reaching $0.325$
at $\alpha=1.9$. The increase is particularly pronounced for smaller
values of $\alpha$, while the CV remains relatively stable between $0.31$
and $0.32$ for $\alpha\geq1.6$. These results indicate that heavier-tailed
L\'evy shocks do not necessarily produce greater relative output
variability. Although smaller $\alpha$ corresponds to heavier tails, the
resulting output dispersion depends on the nonlinear transmission of
stochastic shocks through the saving and capital accumulation dynamics.
Moreover, the CV is estimated from pooled post-burn-in observations,
capturing both within-regime fluctuations and between-regime output
differences. Thus, the observed trend may also reflect changes in capital
regime occupation rather than shock intensity alone. Overall,
Figure~\ref{fig:cv_alpha} demonstrates a nonlinear dependence of
macroeconomic output variability on the L\'evy tail index, highlighting
the importance of endogenous capital dynamics in shaping the response to
stochastic shocks.

\begin{figure}[H]
	\centering
	\includegraphics[width=0.77\textwidth]{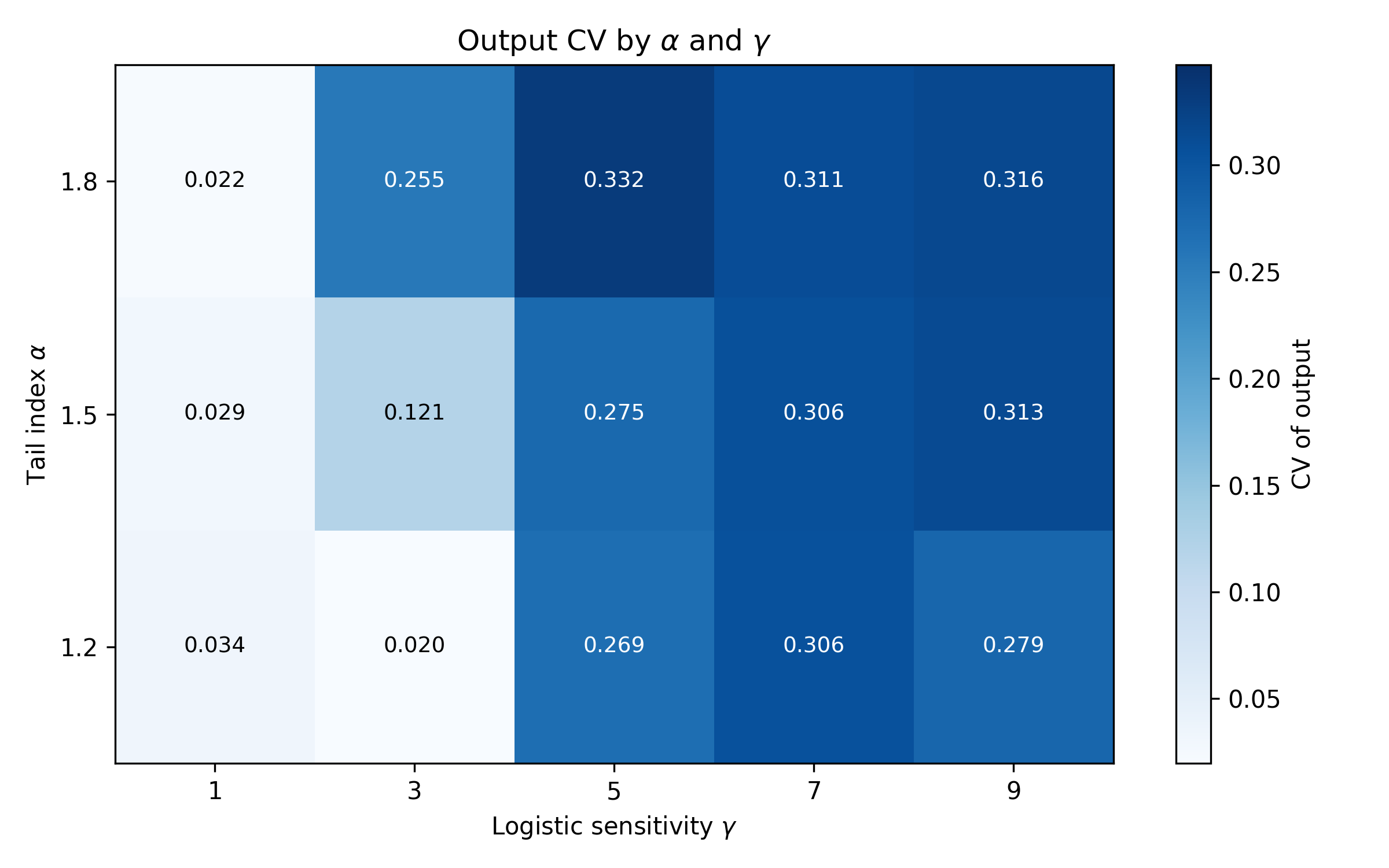}
	\caption{Heat map of the coefficient of variation of output as a function
		of $\gamma$ and $\alpha$.}
	\label{fig:heatmap}
\end{figure}

The joint influence of model parameters and $\alpha$-stable L\'evy
perturbations on stochastic bifurcations and tipping phenomena has been
investigated in nonlinear dynamical systems \cite{yuan2022tipping}. The
interaction between $\gamma$ and $\alpha$ is examined in
Figure~\ref{fig:heatmap}, which displays the estimated coefficient of
variation of output across different parameter combinations. The results
reveal a nonlinear dependence of output variability on both parameters.

Figure~\ref{fig:heatmap} examines the joint effects of the L\'evy
stability index $\alpha$ and the logistic sensitivity $\gamma$ on the CV
of output. The results show that output variability depends nonlinearly
on both parameters. For $\gamma=1$ the CV remains low, ranging from
$0.022$ to $0.034$. As $\gamma$ increases, the CV generally rises,
reaching relatively high values for $\gamma\geq5$. The maximum CV of
approximately $0.332$ occurs at $\alpha=1.8$ and $\gamma=5$, while the
minimum value of $0.020$ occurs at $\alpha=1.2$ and $\gamma=3$. The effect
of $\alpha$ also depends on $\gamma$. For example, at $\gamma=3$ the CV
increases from $0.020$ to $0.255$ as $\alpha$ increases from $1.2$ to
$1.8$. In contrast, at $\gamma=1$ the CV decreases slightly over the same
range. For $\gamma=7$ the CV remains close to $0.31$ for all three values
of $\alpha$. These findings indicate that the effect of L\'evy tail
behavior on output variability depends on the strength of the nonlinear
saving feedback. In particular, heavier-tailed shocks do not necessarily
produce higher output CV. The observed joint sensitivity suggests that the
interaction between stochastic disturbances and endogenous saving dynamics
plays an important role in determining macroeconomic output dispersion.

\subsection{Discussion}\label{sec:discussion}

The numerical results demonstrate that stochastic disturbances to the
scaled saving target are transmitted to the scaled saving rate and
subsequently to capital accumulation through the nonlinear coupled system.
The long-run Monte Carlo distributions exhibit two main concentration
regions associated with the low- and high-capital equilibria. The empirical
occupation fractions and first-passage statistics further indicate that
capital movements depend on both the initial conditions and the L\'evy
stability index. Under the baseline settings, smaller values of $\alpha$
are associated with greater empirical dispersion of the scaled saving
target and more frequent low-to-high first-passage events. However, the
output coefficient of variation exhibits a different pattern. For fixed
$\gamma=5$, the estimated output CV generally increases with $\alpha$ over
the tested range. The joint sensitivity analysis further shows that this
relationship depends on the logistic sensitivity parameter $\gamma$.

These results indicate that heavier-tailed disturbances do not necessarily
produce greater relative output variability. The observed output
dispersion reflects both fluctuations within capital regimes and
differences between the regimes. The occupation probabilities,
first-passage statistics, and output CV are finite-horizon numerical
findings under the specified model parameters. They do not establish
universal monotonicity or global stability properties.

\section{Conclusion and Future Challenges}\label{sec:conclusion}

This work has developed a stochastic growth model in which the saving
rate is an endogenous state-dependent control variable that evolves on a
faster time scale than the capital stock. The model combines three strands
of the growth literature: the Solow tradition of capital accumulation, the
endogenous-saving tradition of Kaldor and Buiter, and the L\'evy-driven
stochastic growth tradition of Abebe et al. \cite{abebe2026}. The central
innovation is the introduction of a logistic saving function that responds
gradually to deviations of capital from a reference level, with the saving
environment subject to heavy-tailed discontinuous shocks.

Section~\ref{sec:stochastic} analyzed the discrete reflected model. We
established invariance of the economically meaningful state space,
derived the stationary distribution of the scaled saving target in the
Gaussian limiting case, and recorded the generator of the
L\'evy-driven system, including the correct factor $h(p)$ in the
compensation term. For $\alpha<2$ no simple closed form for the stationary
distribution is available; the dispersion of the reflected scaled saving
target and the finite-horizon transition statistics are therefore studied
numerically. In particular, we do not claim a general monotonicity theorem
for the stationary variance as a function of $\alpha$; the monotone
decrease observed in the fixed-capital experiments is a numerical finding
for the specified reflected discrete scheme and parameter values.

Section~\ref{sec:numerical} presented numerical simulations that
illustrate the theoretical results and explore the model's implications for
macroeconomic volatility. The simulations confirmed bimodal joint
distributions of capital, scaled saving rate, and scaled saving target, and
asymmetric first-passage behavior between the low- and high-capital
regimes. They also revealed a nonlinear relationship between output
volatility, the logistic sensitivity $\gamma$, and the tail index $\alpha$.
Smaller $\alpha$, corresponding to heavier tails, increased the dispersion
of the scaled saving target and the frequency of low-to-high first-passage
events. However, heavier tails did not uniformly increase relative output
volatility: for fixed $\gamma=5$, the estimated output CV generally
increased with $\alpha$ over the tested range, and the joint heat map
showed that the effect of $\alpha$ depends on $\gamma$. Thus, heavier-tailed
shocks do not necessarily translate into greater relative output
variability; the transmission through state-dependent saving and capital
accumulation is essential.

Several extensions and directions for future research are worth noting.
First, the model assumes that the saving function is logistic, but other
functional forms are possible. For example, a Gompertz or Hill function
could be used to capture different patterns of state-dependent saving.
Second, the model assumes that the noise process is symmetric
$\alpha$-stable, but asymmetric or tempered stable processes could be used
to capture skewness in the shocks. Third, the model assumes that the saving
control and target evolve on a faster time scale than capital, but the
time-scale separation could be relaxed to allow for feedback from capital
to saving on the same time scale. Fourth, the model could be extended to
include endogenous labor supply, technological progress, or open-economy
considerations. Fifth, a rigorous construction of the reflected L\'evy SDE
and a scale-function theory for the coupled three-dimensional system, or a
systematic adiabatic reduction, would provide a useful complement to the
numerical analysis of this paper. Finally, the model could be estimated
using data on saving rates and capital accumulation to assess the
empirical relevance of the state-dependent saving mechanism.

In summary, this work has shown that the interaction between
state-dependent saving and heavy-tailed shocks produces a rich set of
dynamics that are absent in the constant-saving Solow model. The model
provides a unified framework for understanding the joint determination of
saving, capital accumulation, and macroeconomic volatility, and it offers
new insights into the persistence of poverty traps and the sources of
cross-country differences in long-run growth.

\bigskip

\noindent\textbf{Data Availability.}
The numerical algorithms and source code that support the findings of this
study are available from the corresponding author upon reasonable request.

\bigskip

\noindent\textbf{Declaration of competing interest.}
No author associated with this paper has disclosed any potential or
pertinent conflicts which may be perceived as conflicting with this work.

\bigskip

\noindent\textbf{Acknowledgements.}
This work is dedicated to Professor James Brannan on the occasion of his
80th birthday. This work was supported by
Guangdong Basic and Applied Basic Research Foundation (Grant No.
2025A1515012560), Guangdong Introduction Program (Grant No. 2023QN10X753),
and National Foreign Experts Program (Grant No. 111001819820258003).

\end{document}